\documentclass[a4paper]{amsart}

\usepackage{graphicx}
\usepackage{url}

\newcommand{\Dt}{{\Delta t}}
\newcommand{\abs}[1]{\left\vert#1\right\vert}
\newcommand{\R}{\mathbb R}
\newcommand{\N}{\mathbb N}
\newcommand{\Oh}{\mathcal O}
\newcommand{\norm}[1]{\left\Vert#1\right\Vert}
\newcommand{\dott}{\, \cdot\,}
\newcommand{\seq}[1]{\left\{#1\right\}}

\newtheorem{theorem}{Theorem}[section]

\newtheorem{lemma}[theorem]{Lemma}

\newtheorem{remark}[theorem]{Remark}

\usepackage{hyperref}

\allowdisplaybreaks

\begin{document}
\title{Operator splitting for the KdV equation}
\author[Holden]{Helge Holden}
\address[Holden]{\newline
    Department of Mathematical Sciences,
    Norwegian University of Science and Technology,
    NO--7491 Trondheim, Norway,\newline
{\rm and} \newline
  Centre of Mathematics for Applications, 
University of Oslo,
  P.O.\ Box 1053, Blindern,
  NO--0316 Oslo, Norway }
\email[]{holden@math.ntnu.no}
\urladdr{www.math.ntnu.no/\~{}holden}

\author[Karlsen]{Kenneth H.~Karlsen}
\address[Karlsen]{\newline
   Centre of Mathematics for Applications, 
 University of Oslo,
  P.O.\ Box 1053, Blindern,
  NO--0316 Oslo, Norway}
\email[]{kennethk@math.uio.no}
\urladdr{www.math.uio.no/\~{}kennethk}

\author[Risebro]{Nils Henrik Risebro}
\address[Risebro]{\newline
  Centre of Mathematics for Applications, 
University of Oslo,
  P.O.\ Box 1053, Blindern,
  NO--0316 Oslo, Norway}
\email[]{nilshr@math.uio.no}
\urladdr{www.math.uio.no/\~{}nilshr}

\author[Tao]{Terence Tao}
\address[Tao]{\newline
Department of Mathematics, University of California at
  Los Angeles, Los Angeles, CA 90095-1555, USA}
\email[]{tao@math.ucla.edu}
\urladdr{www.math.ucla.edu/\~{}tao}
\date{\today}



\thanks{Supported in part by the Research Council of Norway. 
This paper was written as part of  the international research 
program on Nonlinear Partial Differential Equations at 
the Centre for Advanced Study at the Norwegian Academy 
of Science and Letters in Oslo during the academic year 2008--09. The fourth author is 
supported by a grant from the MacArthur Foundation, the NSF Waterman award, and NSF
grant DMS-0649473.}

\begin{abstract} 
We provide a new analytical approach to operator splitting for equations of the type 
$u_t=Au+B(u)$ where $A$ is a linear  operator and $B$ is  quadratic. 
A particular example is the Korteweg--de Vries (KdV) equation $u_t-u u_x+u_{xxx}=0$. 
We show that the Godunov and Strang splitting methods converge 
with the expected rates if the initial data are sufficiently regular.  
\end{abstract}

\maketitle
\section{Introduction}
The ubiquitous Korteweg--de Vries (KdV) equation 
\begin{equation*}
        u_t-u u_x+u_{xxx}=0
\end{equation*}
offers the perfect blend of the simplest nonlinear convective term $u
u_x$ and the simplest dispersive term $u_{xxx}$: The well-known smooth
soliton solutions of the KdV equation interact in an almost linear
fashion apart from a phase shift.  Furthermore, the KdV equation is
completely integrable with an infinite family of conserved quantities.
This is a result of a subtle interaction between the Burgers term $u
u_x$ and the Airy term $u_{xxx}$, as it is a well-known fact that the
nonlinear Burgers equation $ u_t-u u_x=0$ generically develops shocks
in finite time while the linear Airy equation $u_t+u_{xxx}=0$
preserves all Sobolev norms.

The initial value problem for the KdV equation with $u|_{t=0}=u_0\in
H^s$, either on the whole real line or in the periodic case, has been
extensively studied. An incomplete list of references is
\cite{BonaSmith:1975,KenigPonceVega:1991}. For further discussions we
refer to \cite{Tao:2006,LinaresPonce:2009}.

The method of operator splitting, also called the fractional steps
method, remains a very popular method both for analysis and numerical
computations of partial differential equations. Instead of including a
long list of references to relevant works, we entrust the reader
instead with \cite{HoldenKarlsenLieRisebro:2009} and the overview of
the field given therein.  However, we refer to \cite{Lubich:2008} for
rigorous analysis of the splitting method applied to the
Schr\"odinger--Poisson and the cubic nonlinear Schr\"odinger
equations.

Formally operator splitting can be explained as follows. Let $u(t)=
\Phi_C(t)u_0=\Phi_C(t;u_0)\in X$, where $X$ is some normed space, denote the solution 
of a given differential equation  
\begin{equation*}
  u_t=C(u),\quad t\in[0,T], \qquad u|_{t=0}=u_0. 
\end{equation*}
Here $C$ will typically be a differential operator in the spatial variable. 
Assuming that we can write $C=A+B$ in a natural way, the idea of operator splitting is that
\begin{equation}\label{eq:originalo1}
        u(t_n)\approx \big(\Phi_B(\Dt)\Phi_A(\Dt)\big)^n u_0,
        \quad t_n=n\Dt,
\end{equation}
where $\Dt\ll 1$. In the case of linear ordinary differential
equations this goes back to Sophus Lie. One of the reasons for its
popularity is that operator splitting allows for a separate treatment
of the equations $u_t=A(u)$ and $u_t=B(u)$; in particular this applies
to the use of dedicated special numerical techniques for each of the
equations, again we refer to \cite{HoldenKarlsenLieRisebro:2009} for a
long list of examples and relevant references.

In the context of the KdV equation, the first use of operator
splitting was reported in the brief paper by Tappert
\cite{Tappert:1974} where it was applied as a numerical method.
Apparently the first rigorous results appeared in
\cite{HoldenKarlsenRisebro:1999} where a Lax--Wendroff result was
proved: If operator splitting converges to some limit function, then
the limit function is a weak solution of the KdV equation.  In
addition a systematic study of operator splitting as a numerical
method was undertaken for the KdV equation.  More extensive rigorous
results, specifically the convergence of the splitting approximations
in a suitable functional space, were hampered by the apparent
incompatibility between the Burgers equation and the Airy equation.

In the present paper we offer a new analytical approach to operator
splitting for the KdV equation that will lead to rigorous convergence
results (error estimates). Compared to earlier attempts, two new
ingredients enter the present approach. First of all we actively use
that the solution of the KdV equation remains bounded in a Sobolev
space, that is, if $u_0\in H^s(\R)$, then $\norm{u(t)}_{H^s(\R)}$
remains bounded for $t\in[0,T]$.  This together with a bootstrap
argument is used to secure the existence of a uniform choice of time
step $\Dt$ that prevents the solution from any ``Burgers" step from
blowing up.  Indeed the main problem in this approach is that the Airy
equation produces small oscillatory waves that, when used as initial
data for the Burgers equation, produce shocks. Secondly, since the
splitting approximations are merely defined at the discrete times
$t=t_n$, to facilitate the convergence analysis we introduce an
extension which is defined for all $t\in [0,T]$. Concretely, we
introduce an extension $v$ which depends on an additional time variable
$\tau \in [0,T]$, i.e., $v=v(t,\tau)$, and let the evolution
corresponding to each time variable be governed by one of the split
operators in such a way that at each time level $t=t_n$ the extension
$v(t_n,t_n)$ coincides with the regular splitting approximation.  This
extension approach is 
different from the conventional one, where one lets ``time run twice
as fast'' in each of the subintervals $[t_n,t_n+\Dt/2]$ and
$[t_n+\Dt/2,t_{n+1}]$ (cf.~the discussion and references in
\cite{{HoldenKarlsenLieRisebro:2009}}).

 Formally  the operator splitting \eqref{eq:originalo1}, 
called sequential or Godunov splitting, yields a 
first order approximation in $\Dt$, that is, 
\begin{equation*}
  \norm{u(t_n)-\big(\Phi_A(\Dt)\Phi_B(\Dt)\big)^n u_0}
  =\Oh(\Dt), \quad t_n\to t, \  \Dt\to 0,
\end{equation*}
in an appropriate norm. We show that this holds rigorously for the KdV
equation. More precisely, we prove in Theorem
\ref{thm:godunovconvergence} that for $u_0\in H^s(\R)$ with $s\ge 5$
we have that for $\Dt$ sufficiently small
  \begin{equation*}
    \norm{v(t,t)-u(t)}_{H^{s-3}(\R)}\le K\Dt, \quad t\in [0,T],
  \end{equation*}
where $K$ depends on $s$, $T$ and $u_0$ only. 
Here $v(t,t)$ denotes the splitting approximation in our 
approach evaluated on the ``diagonal" $(t,\tau)=(t,t)$.  

To obtain second order convergence it is 
common to apply the Strang splitting formula, thus formally
\begin{equation*}
        \norm{u(t_n)-\Big(\big(\Phi_B(\Dt/2)\Phi_A(\Dt/2)\big) 
        \big(\Phi_A(\Dt/2)\Phi_B(\Dt/2)\big)\Big)^nu_0}=\Oh(\Dt^2),
\end{equation*}
where $t_n\to t$ as $\Dt\to 0$. 
Here we show rigorously this result for the KdV equation. Indeed, in Theorem 
\ref{thm:Strangsplit} we prove that if $u_0\in H^s$ 
for some $s\ge 17$, then for $\Dt$ sufficiently small 
$$
\norm{v(t,t)-u(t)}_{H^{s-9}} \le K \Dt^2, 
\quad t\in [0,T],
$$
where the constant $K$ depends on $u_0$, $s$ and $T$ only.  Again
$v(t,t)$ denotes the operator splitting approximation in our approach.
Observe that we have to increase the regularity of the initial data,
and hence of the solution, in order to get increased accuracy of
operator splitting.  Note also that with this type of operator
splitting our error estimates are in a much weaker norm than the
assumptions.

It is clear that the present approach applies to several other
equations, and this is currently being investigated.  Furthermore, for
applications to numerical analysis, one would need to replace the
exact solution operators $\Phi_A$ and $\Phi_B$ by numerical
approximations, say $\Phi_A^\delta$ and $\Phi_B^\delta$, and study
their behavior in the limit as $\delta\to 0$, and also to replace the
time derivatives by  discrete differences.
Again this will be studied separately.
 
The paper is organized as follows. We start by presenting the Godunov
operator splitting method for abstract operators. Next we apply this
approach to nonlinear ordinary differential equations where the
procedure is fairly transparent, before we discuss Godunov splitting
for the KdV equation.  Subsequently we present the Strang splitting
for abstract operators. Finally we apply this procedure to the KdV
equation.

\section{Operator splitting}
\label{sec:splitting}
We first present a formal calculation motivating the rigorous analysis
that will follow. Consider an abstract differential equation
\begin{equation}
  \label{eq:original}
  u_t=C(u),\quad t\in[0,T], \qquad u|_{t=0}=u_0 
\end{equation}
for some fixed positive time $T$, 
where $C$ typically will be a differential operator in the spatial
variable. We assume that $u_0$ and $u$ are in some Hilbert space $X$,
and write the solution as
\begin{equation*}
  u(t)=\Phi_C(t;u_0).
\end{equation*}
Formally, expanding the solution in a Taylor series we find
\begin{equation*}
  u(t)=u_0+tu_t(0)+\Oh(t^2)= u_0+t C(u_0)+\Oh(t^2),
\end{equation*}
where we have used the equation \eqref{eq:original}. 
Assume that one can write
\begin{equation*}
  C=A+B
\end{equation*}
in some natural way. Operator splitting (of the Godunov type)
then works as follows: Instead of solving the
problem with $C$ directly, one alternately solves for small time steps
the equations $u_t=A(u)$ and $u_t=B(u)$. 
Making the time steps finer and finer, the approximation 
will presumably converge to the solution of the original equation \eqref{eq:original}.

More precisely, fix a positive time step $\Dt$, let $t_n=n \Dt$, 
$n\in\N_0$, and define a family 
$\{u_\Dt(t_n)\}_n$ of functions
\begin{equation*}
  \begin{aligned}
    u_{\Dt}(t_{n+1})&= \Phi_A(\Dt;\Phi_B(\Dt;u_{\Dt}(t_{n}) ))=
    \Phi_A(\Dt)\circ\Phi_B(\Dt)u_{\Dt}(t_{n}), \quad n\in\N_0, \\ 
    u_{\Dt}(0)&=u_0.
  \end{aligned}
\end{equation*}
The traditional method of extending the solution to any $t\in [0,T]$
has been to let ``time run twice as fast'' in each of the subintervals
$[t_n,t_{n+1/2}]$ and $[t_{n+1/2},t_{n+1}]$, 
where $t_{n+1/2}=t_n+\Dt/2$. Thus obtaining
\begin{equation}
  u_{\Dt}(t)=
  \begin{cases}
    \Phi_{B}\left(2(t-t_n);u_{\Dt}\left(t_n\right)\right)
    &\text{for $t\in \left[t_n,t_{n+1/2}\right]$,}\\
    \Phi_{A}\left(2(t-t_{n+1/2});u_{\Dt}(t_{n+1/2})\right)
    &\text{for $t\in \left[t_{n+1/2},t_{n+1}\right]$.}
  \end{cases}\label{eq:utraddef}
\end{equation}
(In the present approach we will use a different extension to all times $t$.)

Formally one can show
\begin{equation*} 
        \norm{u_{\Dt}(t_n)-u(t_n)}\le \Oh(\Dt) 
        \text{ as $\Dt\to 0$ and $t_n\to t$}
\end{equation*}
in some norm. 

The convergence can be improved to second order by using the 
Strang splitting formula. To this end we let the approximation (this
time denoted $v_\Dt$ to distinguish from the previous approximation) read
\begin{equation*}
  \begin{aligned}
    v_{\Dt}(t_{n+1})&=\Phi_B(\Dt/2;\Phi_A(\Dt;\Phi_B(\Dt/2;v_{\Dt}(t_{n}))))\\
    &= \Phi_B(\Dt/2)\circ \Phi_A(\Dt)\circ \Phi_B(\Dt/2)
    v_{\Dt}(t_{n})\\
    &=\Big(\Phi_B(\Dt/2)\circ \Phi_A(\Dt/2)\Big)
        \\ & \qquad\quad \circ\Big(
    \Phi_A(\Dt/2)\circ \Phi_B(\Dt/2)\Big)
    v_{\Dt}(t_{n}) , \quad n\in\N_0, \\
    v_{\Dt}(0)&=u_0.
  \end{aligned}
\end{equation*}
Formally one now has
\begin{equation*}
  \norm{v_{\Dt}(t_{n})-u(t_n)}\le \Oh(\Dt^2) \text{ as $\Dt\to 0$ and $t_n\to t$.}
\end{equation*}

To show that the operator splitting solutions are well defined, we shall later 
make use of the following bootstrap lemma,
taken from \cite[Prop.~1.21]{Tao:2006}. 
\begin{lemma}
  \label{lem:bootstrap} Let $t\in [0,T]$. Consider a continuous function $\phi\colon [0,T] \to[0,\infty)$. 
If there exists a positive constant $\alpha$ such that:
  \begin{itemize}
  \item[({\bf a})] $\phi(0)\le \alpha$;

  \item[({\bf b})] for any $t$ such that $\phi(t) \le \alpha$, we can show that
     $\phi(t)\le \alpha/2$; 
  \end{itemize}
  then $\phi(t)\le \alpha/2$ for all $t\in [0,T]$.
\end{lemma}

\subsection{Doubling the time variable}\label{subsec:tdouble}
We shall formulate the operator splitting solution by introducing two time variables, and 
define a function $v= v(t,\tau)$ for $(t,\tau)$ in the set
\begin{equation*}
        \Omega_{\Dt}= \bigcup_{n=0}^{\lfloor T/\Dt\rfloor}
        [t_n,t_{n+1}]\times [t_n,t_{n+1}]
\end{equation*}
by requiring  that
\begin{equation}
  \label{eq:v}
  \begin{aligned}
    v(0,0)&=u_0,\\
    v_t(t, t_n)&=B(v(t,t_n)), \quad t\in (t_n,t_{n+1}], \\
    v_\tau(t,\tau)&=A(v(t,\tau)), \quad (t,\tau)\in [t_n,t_{n+1}]\times(t_n,t_{n+1}],
\end{aligned}
\end{equation}
where $n=0,\dots,\lfloor T/\Dt\rfloor$.
\begin{figure}[tbp]
  \centering
  \includegraphics[width=0.5\linewidth]{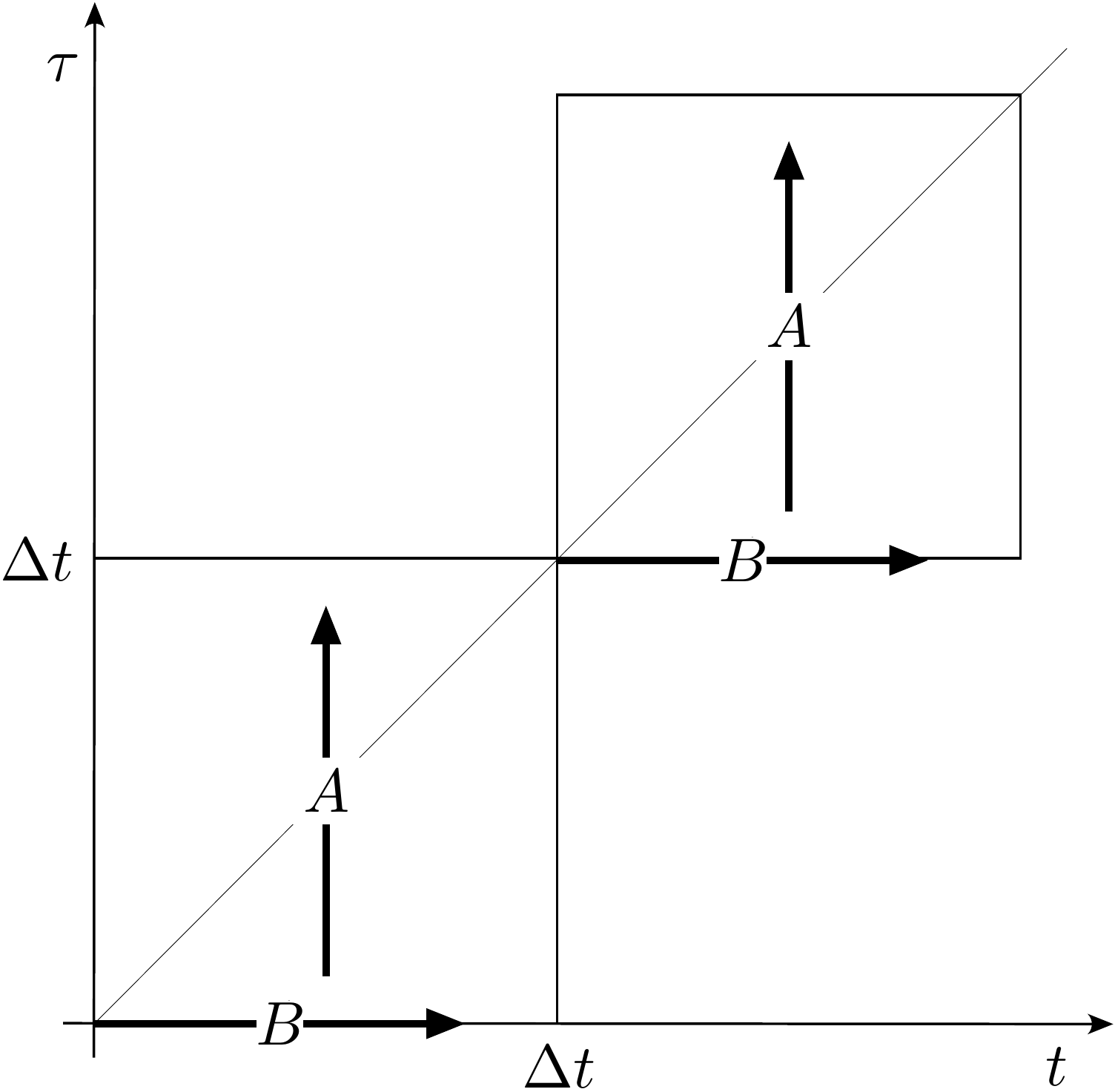}
  \caption{Schematic view of the Godunov splitting and the definition
    of $v(t,\tau)$, cf.~equation \eqref{eq:v}.} 
  \label{fig:godsplit}
\end{figure}
Observe that 
\begin{equation*}
  u_{\Dt}(t_n)=v(t_n,t_n),\quad n=0,\dots, \lfloor T/\Dt\rfloor. 
\end{equation*}
where $u_{\Dt}$ is given by \eqref{eq:utraddef}. 
The specific extension $v$ of $\{u_\Dt(t_n)\}_n$ to $[0,T]$ 
will serve as an important technical tool to be utilized in the analysis.

The exact solution of \eqref{eq:original} is still denoted by $u$. 
Introduce the error function
\begin{equation*}
  w(t)=v(t,t)-u(t).  
\end{equation*}
The aim is to show that 
\begin{equation*}
  \norm{w(t)}\le\Oh(\Dt), \quad  \Dt\to 0,\quad t\in[0,T]
\end{equation*}
in an appropriate norm.

We introduce second order Taylor expansions of the operators 
$A$ and $B$ (see~\cite[p.~29]{AmbrosettiProdi}),
\begin{align*}
  A(f+g)&=A(f)+dA(f)[g]+\int_0^1(1-\alpha)
  d^{(2)}A(f+\alpha g)[g]^2 d\alpha, \\
  B(f+g)&=B(f)+dB(f)[g]+\int_0^1(1-\alpha)
  d^{(2)}B(f+\alpha g)[g]^2 d\alpha.
\end{align*}
Thus we find
\begin{equation}
  \label{eq:w_eq}
  \begin{aligned}
    w_t- dA(u)[w]-dB(u)[w]&= v_t+v_\tau-u_t-dA(u)[w]-dB(u)[w] \\
    &=v_t+A(v)-(A+B)(u)- dA(u)[w]-dB(u)[w] \\
    &=v_t-B(v)+\big(A(v)-A(u)-dA(u)[w]\big)\\
    &\qquad +\big(B(v)-B(u)-dB(u)[w]\big)\\
    &=F(t)+ \int_0^1 (1-\alpha) d^{(2)}A(u+\alpha w)[w]^2 d\alpha\\
    &\qquad +\int_0^1(1-\alpha) d^{(2)}B(u+\alpha w)[w]^2 d\alpha
  \end{aligned}
\end{equation}
where we have introduced a forcing term $F(t)=F(t,t)$ defined as
\begin{equation}
  \label{eq:forc}
  F(t,\tau)=v_t(t,\tau)-B(v(t,\tau)).
\end{equation}
We can rewrite \eqref{eq:w_eq} as follows
\begin{equation}
  \label{eq:w_eqA}
 w_t- dC(u)[w]= F(t)+ \int_0^1 (1-\alpha) d^{(2)}C(u+\alpha w)[w]^2 d\alpha.
\end{equation}
The forcing term satisfies the following time development
\begin{equation}
  \label{eq:F_tau}
  \begin{aligned}
    F_\tau-dA(v)[F]&=v_{t\tau}-B(v)_\tau-dA(v)[v_t-B(v)]\\
    &=A(v)_t-dB(v)[v_\tau]-dA(v)[v_t]+dA(v)[B(v)]\\
    &=dA(v)[v_t]- dB(v)[A(v)]-dA(v)[v_t]+dA(v)[B(v)] \\
    &=[A,B](v,v),
  \end{aligned}
\end{equation}
where we have defined the commutator
\begin{equation*}
  [A,B](f,g)=dA(f)[B(g)]-dB(f)[A(g)].
\end{equation*}
For simplicity we will subsequently be writing $[A,B](v)$ rather than $[A,B](v,v)$.

\subsection{Ordinary differential equations} \label{subsec:odes} 
As a warm-up we consider the ordinary differential equation
\begin{equation}
u_t=C(u), \quad t>0, \quad u(0)=u_0\in \R^n.\label{eq:odeex}
\end{equation}

To simplify the presentation we assume that $C$ is 
a quadratic function, i.e., that $C'''=d^3C=0$ or $d^2C[f,g]$ is constant. 
This means that the integral in \eqref{eq:w_eqA} reduces 
to a constant. Furthermore, we assume that  \eqref{eq:odeex} is such that
there is a unique solution $u(t)$ such that $\abs{u(t)}\le K_{u_0,T}$
for $t\in [0,T]$. Furthermore, we will assume 
that the operators $A$ and $B$ are two 
times continuously differentiable and
\begin{equation}
  \label{eq:assu_ODE}
   A',A''\in L^\infty,\quad A(0)=0, \qquad \abs{B(u)}\le K \abs{u}^2, \quad \abs{B'(u)}\le K \abs{u}.
\end{equation}

Throughout this paper we use the convention that for a quantity
$\alpha$, $K_\alpha$ denotes a constant depending on
$\alpha$ (and perhaps other things). We use this notation to highlight the
dependence on $\alpha$. The actual value of $K_\alpha$ may be different at each
occurrence.

Let now $\alpha$ be a positive constant (its precise value will be fixed
later). To start the bootstrap argument we assume
\begin{equation}\label{eq:Aodeass}
        \abs{v(t,\tau)}\le \alpha, \quad (t,\tau)\in\Omega_\Dt.
\end{equation}
In this example, since $C$ is quadratic, $d^2C=\kappa$ for some
constant symmetric matrix $\kappa$, and $w$ satisfies 
\begin{equation}
        w_t + C'(u)w = F + w^T\frac{\kappa}{2} w,
        \quad t>0,\quad w(0)=0.\label{eq:wodedev}
\end{equation}
Furthermore, $F$ satisfies 
$$
F_\tau - A'(v)F = [A,B](v), \quad 
(t,\tau)\in [t_n,t_{n+1}]\times (t_n,t_{n+1}],
$$
and $F(t,t_n)=0$, for each $n$.
This means that 
$$
\frac{\partial}{\partial \tau} \abs{F}
\le K \abs{F} + K\alpha^2.
$$
Hence, Gronwall's inequality implies that 
$$
\abs{F}(t,\tau) \le K_\alpha \Dt,
$$
where $K_\alpha$ is a constant depending on the assumed 
bound on $v$ in \eqref{eq:Aodeass}. 
In view of this bound and \eqref{eq:wodedev},
$$
\frac{d}{dt}\abs{w} \le K_\alpha \abs{w} + K_\alpha\Dt, \quad w(0)=0.
$$
Gronwall's inequality gives that 
$$
\abs{w(t)} \le e^{K_\alpha t} t K_\alpha \Dt \le K_\alpha \Dt.
$$

Trivially we have
$$
\abs{v(t,t)-v(t,\tau)}\le 
\int_{\min\seq{t,\tau}}^{\max\seq{t,\tau}}
\abs{A(v(t,s))}\,ds \le K_\alpha \Dt,
$$
for $(t,\tau)\in [t_n,t_{n+1}]\times (t_n,t_{n+1}]$ for any $n$.
Then we can conclude that
$$
\abs{v(t,\tau)} \le \abs{u(t)} 
+ \abs{w(t)}+\abs{v(t,t)-v(t,\tau)}  
\le K + K_\alpha \Dt.
$$
Now we are in a position to choose $\alpha$ 
so that $K\le \alpha/4$, this determines $K_\alpha$. 
Next choose $\Dt$ so small that $K_\alpha\Dt\le \alpha/4$, then
$$
\abs{v(t,\tau)}\le \alpha/2  \quad (t,\tau)\in\Omega_\Dt.
$$
Hence, by the bootstrap lemma, $\abs{v(t,t)}\le \alpha/2$ for all
$t$. Consequently, 
\begin{equation} \label{eq:ODE_G}
        \abs{v(t,t)-u(t)}\le K_{\alpha/2} \Dt,
\end{equation}
i.e., the operator splitting is as expected first order accurate.

\begin{remark}\label{rem:ODE1}
An interesting example is the logistic equation $u'=u(u-1)$ where we
can write $A(u)=-u$ and $B(u)=u^2$. Exact solutions are available for
all operators involved, specifically
\begin{equation*}
  \Phi_C(t)u_0=\frac{u_0}{u_0+e^t(1-u_0)}, \quad \Phi_A(t)u_0=u_0e^{-t}, 
\quad \Phi_B(t)u_0=\frac{u_0}{1-u_0t}.
\end{equation*}
Let $u_0\in(0,1)$. Then there is no blow-up in the full equation, 
but blow-up for the equation $u_t=b(u)$ at $t^*=1/u_0$.
The function $v$ reads in our case
\begin{equation*}
        v(t,\tau)=\frac{v(t_n,t_n)e^{-(\tau-t_n)}}{1-v(t_n,t_n)(t-t_n)}, 
        \quad t,\tau\in[t_n,t_{n+1}], 
\end{equation*}
where
\begin{equation*}\label{eq:ex1a_ODE}
        v(t_n,t_n)= \frac{u_0(1-e^{-\Dt})}{(1-e^{-\Dt})e^{t_n}+u_0\Dt(1-e^{t_n})}. 
\end{equation*}
\begin{figure}[tbp]
  \centering \includegraphics[width=0.45\linewidth]{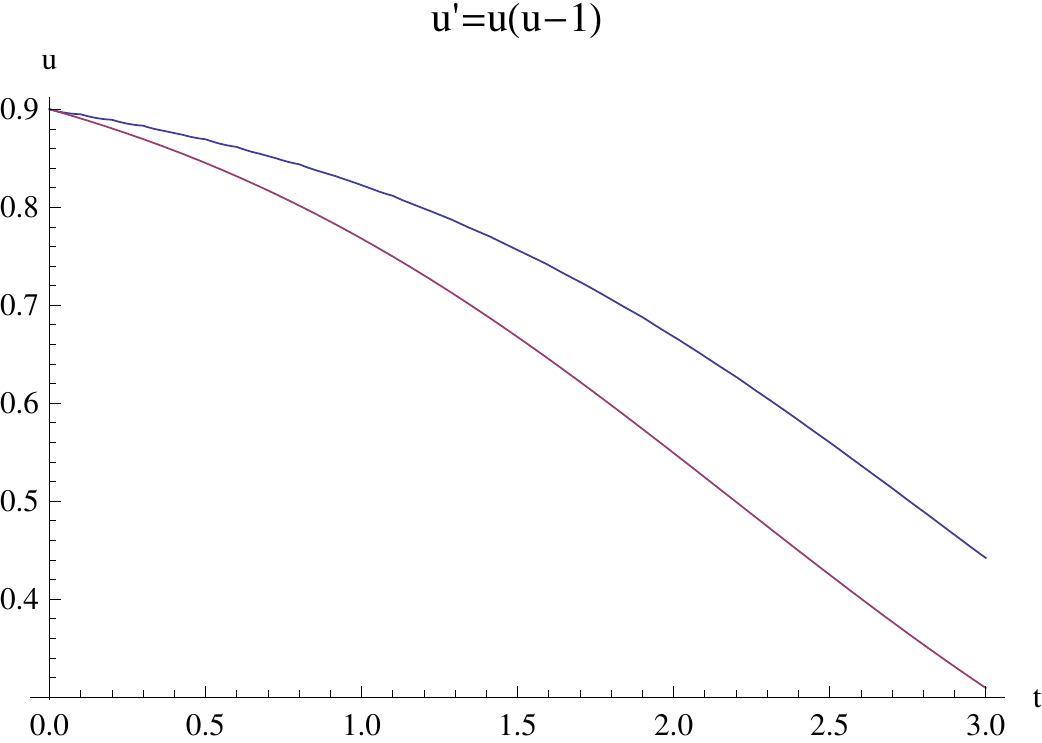}
  \includegraphics[width=0.45\linewidth]{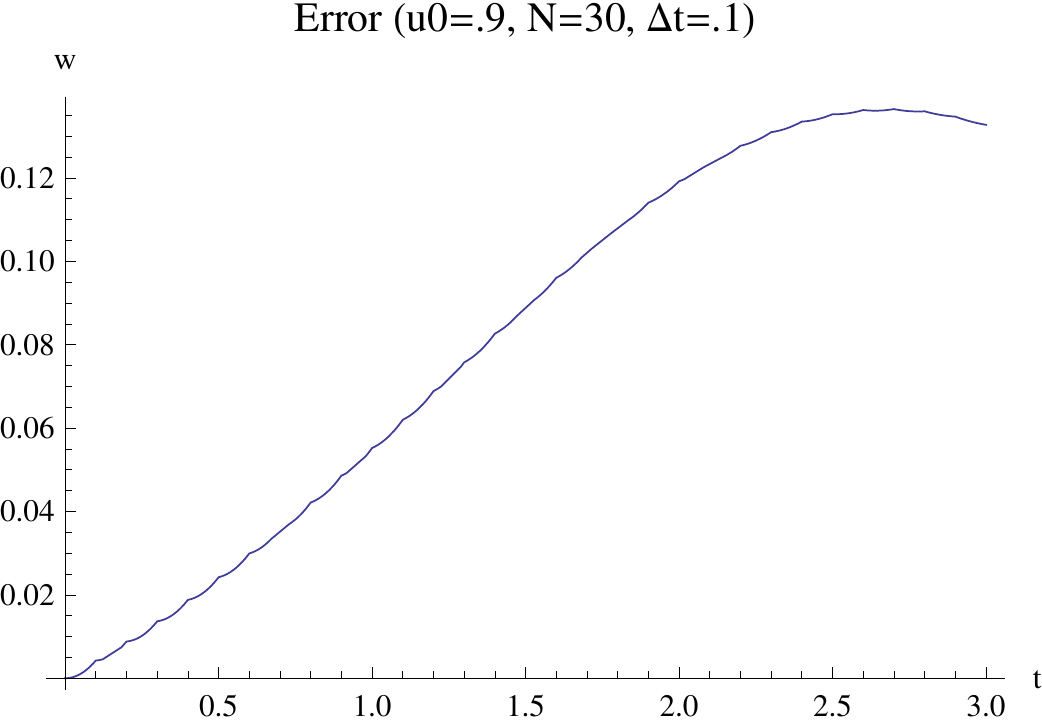}\\
  \includegraphics[width=0.45\linewidth]{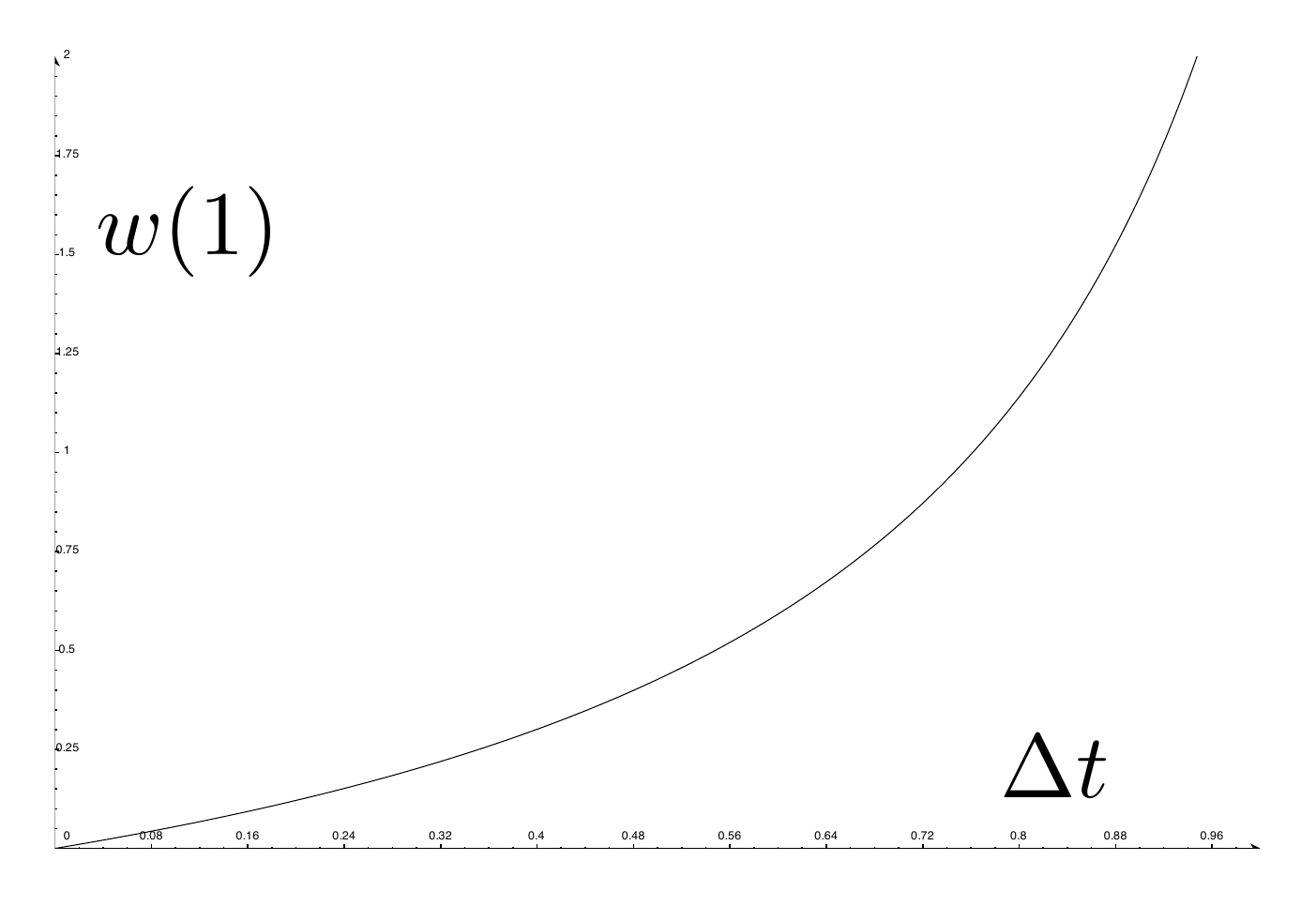}
  \caption{(Left) The exact (blue) and the approximate (red)
    solution. (Right) The error. (Below) The error $w(1)$ as a
    function of $\Dt$.}
  \label{fig:odeSplit}
\end{figure}

The quantity $v(t_n,t_n)$ is well-defined on $[0,T]$ if one 
chooses $\Dt$ such that
\begin{equation*}\label{eq:Dt}
        \Dt < \frac{1-e^{-\Dt}}{u_0(1-e^{-t_n})}.
\end{equation*}
Since $t_n\le T$, and  $e^T < u_0/(u_0-1)$, we find that
$u_0(1-e^{-t_n})\le u_0(1-e^{-T}) <1$. 
Thus we have to choose $\Dt$ such that
\begin{equation*}
        u_0\big(1-e^{-T}\big)<\frac1{\Dt}\big(1-e^{-\Dt}\big)= 1-\frac{\Dt}{2}+\Oh\big(\Dt^2\big), 
\end{equation*}
or 
\begin{equation*}\label{eq:Dt2a}
        \Dt< 2\big(1- u_0(1-e^{-T})\big). 
\end{equation*}
In this case one can verify \eqref{eq:ODE_G} directly, namely
\begin{equation*}
        \begin{aligned}
                & \abs{v(t_n,t_n)-u(t_n)}\\
                &\quad=\abs{u_0}^2(e^{t_n}-1)
                \frac{\abs{(\Dt)^{-1}(1-e^{-\Dt})-1}}{\abs{
                (u_0+e^{t_n}((\Dt)^{-1}(1-e^{-\Dt})-u_0))(u_0+e^{t_n}(1-u_0))}}\\
                & \quad =\Oh(\Dt).
        \end{aligned}
\end{equation*}
The example is illustrated in Figure \ref{fig:odeSplit}.
\end{remark}

\subsection{The KdV equation}
\label{sec:kdv_split}
Let us now apply this general framework to the KdV
equation, that is, 
\begin{equation*} 
  u_t=uu_x-u_{xxx}, \quad u|_{t=0}=u_0 \in H^s(\R).
\end{equation*}
In this case $C(u)=uu_x - u_{xxx}$, and the evolution operator
$\Phi_C(t;\dott)\colon H^s \to H^s$ is bounded. Therefore, the Hilbert space
$\R^n$ of the previous example is replaced by the Sobolev space
$H^s(\R)$ with the inner product and norm
$$
(f,g)_{H^s(\R)} = \sum_{j=0}^s \int_{\R} \partial^j_x f(x)\, \partial^j_x
g(x)\, dx, \quad \norm{f}^2_{H^s(\R)} = (f,f)_{H^s(\R)}.
$$

We choose $A$ to  equal (minus) the Airy operator:
\begin{equation*}
  \begin{aligned}
    A(f)&=-f_{xxx}, \\
    dA(f)[g]&=-g_{xxx}, \\
    d^{(2)}A(f)[g,h]&=0, 
  \end{aligned}
\end{equation*}
and $B$ to equal  the Burgers operator:
\begin{equation*}
  \begin{aligned}
    B(f)&=ff_{x}, \\
    dB(f)[g]&=fg_{x}+f_xg, \\
    d^{(2)}B(f)[g,h]&=hg_{x}+h_xg,\\
    d^{(3)}B(f)[g,h,k]&=0.
  \end{aligned}
\end{equation*}
In this case the commutator reads
\begin{equation*}
  [A,B](f,f)=-\frac32 \partial_x^2(f_x)^2.  
\end{equation*}
Thus the equations \eqref{eq:w_eqA}, \eqref{eq:forc}, and \eqref{eq:F_tau} are
\begin{align}
        w_t-(u w)_x+w_{xxx}&=F+w w_x,
        \label{eq:kdvW} \\
        F&=v_t-vv_x, \notag 
        \\ F_\tau+F_{xxx}&=-\frac32 \partial_x^2(v_x)^2,
        \label{eq:kdvFtau}
\end{align}
respectively. 

From  \cite{BonaSmith:1975} we recall the classical result: For $u_0\in
H^s(\R)$ with $s\ge 2$ there exists a unique solution $u\in C([0,T],
H^s(\R))$ of the KdV equation
\begin{equation*}
  u_t=u u_x-u_{xxx}, \quad u|_{t=0}=u_0.
\end{equation*}
In particular, we can assume that there exists a constant $K$
(depending on $T$, $u_0$, and $s$) such that
\begin{equation*}
  \norm{u(t)}_{H^s(\R)}\le K, \quad t\in [0,T]. 
\end{equation*}
To save space and typing efforts, we shall write $H^k$ for $H^k(\R)$,
and $\partial$ for $\partial_x$.

Next, let $s$ be an odd integer larger or equal to $5$. Assume 
that there exists a constant $\alpha$ such that
\begin{equation*}
  \norm{v(t,\tau)}_{H^{\hat{k}}}\le \alpha, \quad (t,\tau)\in\Omega_\Dt,
\end{equation*}
where $\hat{k}(s)=(s-1)/2$. 

Let us estimate the behavior of the Airy and Burgers
operators. The Airy equation leaves all Sobolev norms invariant, viz.
\begin{equation*}
  \norm{v(t,\tau)}_{H^k}=\norm{v(t,t_n)}_{H^k}. 
\end{equation*}

By definition we find for the Burgers operator
\begin{align*} 
  \frac12\frac{d}{dt} \norm{v(t,t_n)}_{H^s}^2&=(v,v_t)_{H^s} 
  = \sum_{j=0}^s \int \partial^j v\partial^j\left(vv_x\right)\, dx
  \\
  &= \sum_{j=0}^s \sum_{k=0}^j \binom{j}{k} \int
  \partial^j v \, \partial^{k+1}v\, \partial^{j-k}v \,dx.
\end{align*}
For $j<s$, any of the above terms can be estimated by
\begin{align*}
  \Bigl| \int \partial^j v\,\partial^{k+1} v\, \partial^{j-k}
  v\,dx\Bigr|&\le 
  \norm{\partial^{j} v}_{L^\infty}
  \norm{\partial^{\max\seq{k+1,j-k}} v}_{L^2} 
  \norm{\partial^{\min\seq{k+1,j-k}} v}_{L^2}
  \\
  &\le K \norm{v}_{H^s}^2 \norm{v}_{H^{\hat{k}}} \\
  &\le K_\alpha \norm{v}_{H^s}^2.
\end{align*}
For $j=s$ all terms with $k<s$ can be estimated similarly,
\begin{align*}
  \Bigl| \int \partial^s v\,\partial^{k+1} v\, \partial^{s-k}
  v\,dx\Bigr|&\le 
  \norm{\partial^{s} v}_{L^2}
  \norm{\partial^{\min\seq{k+1,s-k}} v}_{L^\infty}
  \norm{\partial^{\max\seq{k+1,s-k}} v}_{L^2} 
  \\
  &\le K \norm{v}_{H^s}^2 \norm{v}_{H^{\hat{k}}} \\
  &\le K_\alpha \norm{v}_{H^s}^2.
\end{align*}
We
are left with the term where $k=s=j$, viz.
\begin{align*}
  \Bigl| \int \partial^s v\, \partial^{s+1} v \, v \,dx \Bigr| &=
  \frac{1}{2}\Bigl | \int \left(\partial^s v\right)^2 \, \partial
  v\,dx\Bigr| \\
  &\le \norm{\partial v}_{L^\infty} \norm{\partial^s v}_{L^2}^2\\
  &\le K \norm{v}_{H^{\hat{k}}} \norm{v}_{H^s}^2 
\end{align*}
if $\hat{k}\ge 2$, i.e., if $s\ge 5$. 

Thus 
\begin{equation}
\frac{d}{dt} \norm{v\left(t,t_n\right)}_{H^s} \le K_\alpha
\norm{v\left(t,t_n\right)}_{H^s},\label{eq:burgerhsdev}
\end{equation}
which implies that
$$
\norm{v\left(t,t_n\right)}_{H^s} \le e^{K_\alpha(t-t_n)}
\norm{v\left(t_n,t_n\right)}_{H^s}.
$$
In particular, for any $n$,
\begin{equation*}
        \norm{v\left(t_n,t_n\right)}_{H^s}\le e^{K_\alpha\Dt}
        \norm{v\left(t_{n-1},t_{n-1}\right)}_{H^s} 
        \le e^{K_\alpha t_n} \norm{u_0}_{H^s}.
\end{equation*}

Thus we have shown the following result.
\begin{lemma}\label{lem:indu}
We have
\begin{equation*}
        \norm{v(t,\tau)}_{H^s}\le K_\alpha, \quad (t,\tau)\in\Omega_\Dt.
\end{equation*}
\end{lemma}

Observe the general result, obtained by integration by parts,
$$
(f,f_{xxx})_{H^k}=0, \quad f\in H^{k+3}.
$$
Next we analyze the forcing term that satisfying \eqref{eq:kdvFtau}. By
taking the $H^{s-3}$ inner product with $F$ in \eqref{eq:kdvFtau} we get
\begin{align*} 
  \frac12 \partial_\tau \norm{F(t,\tau)}_{H^{s-3}}^2 &=
 -(F,F_{xxx})_{H^{s-3}} -\frac32(F,\partial^2(v_x)^2)_{H^{s-3}}\\
  &\le \frac32 \norm{F}_{H^{s-3}}\norm{\partial_x^{j+2}(v_x)^2}_{H^{s-3}}.
\end{align*}
Since $H^s$ is an algebra
\begin{equation*}
    \norm{\partial^{j+2}(v_x)^2}_{H^{s-3}}\le \norm{(v_x)^2}_{H^{s-1}}
    \le K \norm{v_x}_{H^{s-1}}^2
    \le K_\alpha^2.   
\end{equation*}
Also $F(t,t_n)=0$ for $t\in[t_n,t_{n+1}]$, and we conclude that
\begin{equation*}
  \norm{F(t,\tau)}_{H^{s-3}}\le K_\alpha \Dt, \quad (t,\tau)\in \Omega_\Dt.
\end{equation*}

As for the error function $w=v(t,t)-u(t)$, we have the following estimates.  Let 
$E(t)=\norm{w(t)}_{H^{s-3}}$. By taking the inner product with $w$
in \eqref{eq:kdvW} 
\begin{equation}
  \label{eq:west}
  \begin{aligned}
    \frac{1}{2}\frac{d}{dt} E^2(t)&
=(w,(uw)_x)_{H^{s-3}}-(w,w_{xxx})_{H^{s-3}}+(w,F)_{H^{s-3}}+(w,ww_x)_{H^{s-3}} \\  
&=  (w,(uw)_x)_{H^{s-3}} +(w,ww_x)_{H^{s-3}} +(w,F)_{H^{s-3}}\\
  &  =\sum_{j=0}^{s-3}\int \left(\partial^j w\partial^j\left(u w_x+u_xw+w
        w_x\right) \right)\, dx+(w,F)_{H^{s-3}}.
      \end{aligned}
\end{equation}
The first integrand on the right is expanded by Leibniz' rule. 
We get
$$
\sum_{j=0}^{s-3} \sum_{k=0}^{j} 
\binom{j}{k} \int \partial^j w \partial^{k+1}w \partial^{j-k} u + 
\partial^j w \partial^k w \partial^{j+1-k} u + \partial^j
w \partial^{k} w \partial^{j+1-k} w \,dx.
$$
For $0\le k \le j < s-3$, we can estimate 
\begin{align*}
  \Bigl| \int \partial^j w \, \partial^{k+1 (k)} w \, \partial^{j-k (+1)} 
  u \,dx\Bigr| &\le  \norm{w}_{H^{s-3}}^2
  \norm{\partial^{j-k (+1)}u}_{L^\infty}  \\
&  \le K \norm{w}_{H^{s-3}}^2 \norm{u}_{H^s}, \\[1mm]
  \Bigl| \int \partial^j w \,\partial^k w \, \partial^{j+1-k} w
  \,dx\Bigr| 
  &\le \norm{\partial^j w}_{L^\infty} \norm{w}_{H^{s-3}}^2\\
  &\le
  K\left(\norm{u}_{H^s} + \norm{v}_{H^s}\right)\norm{w}_{H^{s-3}}^2.
\end{align*}
For $j=s-3$, we can use the same strategy for those terms with fewer
than $s-2$ derivatives on $w$. The term with $s-2$ derivatives on $w$
can be estimated as
\begin{align*}
  \Bigl| \int \partial^{s-3}w \, \partial^{s-2} w \,u \,dx \Bigr| &=
  \frac{1}{2} \Bigl|\int \left(\partial^{s-3} w\right)^2 \, \partial u \,dx
  \Bigr| \le K\norm{u}_{H^s} \norm{w}_{H^{s-3}}^2,\\[2mm]
  \Bigl| \int \partial^{s-3} w \, w \partial^{s-2} w \,dx\Bigr| &=
  \frac{1}{2} \Bigl| \int \left(\partial^{s-3} w\right)^2 \,\partial
  w\,dx\Bigr| \\
 & \le K\left(\norm{u}_{H^s}+\norm{v}_{H^s}\right)\norm{w}_{H^{s-3}}^2.
\end{align*}
The last term in \eqref{eq:west} is overestimated by
$\norm{F}_{H^{s-3}} \norm{w}_{H^{s-3}}$, and we get
\begin{equation*}
  \frac{d}{dt}E^2(t) \le K_\alpha E^2(t)+ K E(t) \norm{F}_ {H^{s-3}}
\le  K_\alpha E^2(t)+ K \Dt,
\end{equation*}
which implies
\begin{equation*}
  \frac{d}{dt} E(t) \le K_\alpha E(t)+ K \Dt.
\end{equation*}
Since $E(0)=0$, Gronwall's inequality yields
\begin{equation*}
  \norm{w(t)}_{H^{s-3}}=E(t) \le K_\alpha \Dt.
\end{equation*}
Recall that $\norm{v(t,\tau)}_{H^{s-3}} =\norm{v(t,t)}_{H^{s-3}}$
because
$$
\partial_\tau \norm{v(t,\tau)}_{H^{s-3}}^2 =
-2\left(v,v_{xxx}\right)_{H^{s-3}}=0.
$$
Now $u$ is bounded in $H^s$, and we infer that
\begin{equation*}
  \norm{v(t,\tau)}_{H^{s-3}}\le K+ K_\alpha \Dt. 
\end{equation*}
Since $s\ge 5$,
$$
s-3\ge \hat{k}=\frac{s-1}{2},
$$
and we get
\begin{equation*}
  \norm{v(t,\tau)}_{H^{\hat{k}}}\le K+ K_\alpha \Dt.
\end{equation*}
First choose $\alpha\ge 4K$, then choose $\Dt$ so small that $K_\alpha\Dt \le
\alpha/4$, so that $K+K_\alpha\Dt\le \alpha/4+ \alpha/4= \alpha/2$. Hence, by the bootstrap argument
we have proved the following theorem:
\begin{theorem}
  \label{thm:godunovconvergence}
  Fix $T>0$. Let $u_0\in H^s(\R)$ with $s\ge 5$. Then for $\Dt$
  sufficiently small we have 
  \begin{equation*}
    \norm{v(t,t)-u(t)}_{H^{s-3}(\R)}\le K\Dt, \quad t\in [0,T],
  \end{equation*}
where $K$ depends on $s$, $T$ and $u_0$ only.
\end{theorem}
\begin{remark} \label{rem:interchange} 
  Instead of defining $v$ by \eqref{eq:v} for the KdV equation, we
  could also interchange the order of the Airy operator $A$ and the
  Burgers operator $B$ in the definition of $v$.  The same procedure
  as described above would apply, and Theorem
  \ref{thm:godunovconvergence} would remain valid. This remark is
  important for the Strang splitting to be discussed next.
\end{remark}

\section{Strang splitting}
\label{sec:strang}

To achieve higher-order convergence it is common to consider the so-called
Strang splitting. Now we approximate the solution by using two Godunov
splittings, each with a time step of $\Dt/2$, and in alternating
order. Explicitly, we define
\begin{align}
  v(0,0)&=u_0,\notag\\
  v_t(t, t_n)&=B(v(t,t_n)), \quad t\in (t_n,t_{n+1/2}],
 \notag \\
  v_\tau(t,\tau)&=A(v(t,\tau)), \quad (t,\tau)\in
  [t_{n},t_{n+1/2}]\times (t_n,t_{n+1/2}], \label{eq:GodunovDef}
  \\
  v_\tau(t_{n+1/2},\tau)&=A(v(t_{n+1/2},\tau)), \quad \tau\in
  (t_{n+1/2},t_{n+1}], 
 \notag \\
  v_t(t, \tau)&=B(v(t,\tau)), \quad (t,\tau)\in
  (t_{n+1/2},t_{n+1}]\times [t_{n+1/2},t_{n+1}], \notag
\end{align}
for $ n=0,\dots, \lfloor T/\Dt\rfloor$. We consider this function for
$(t,\tau)$ in the domain
\begin{equation*}
  \widetilde\Omega_{\Dt}= \bigcup_{n=0}^{\lfloor T/\Dt\rfloor}
 \left( \left[t_n,t_{n+1/2}\right]^2 \cup
  \left[t_{n+1/2},t_{n+1}\right]^2\right).
\end{equation*}
\begin{figure}[tbp]
  \centering
  \includegraphics[width=0.5\linewidth]{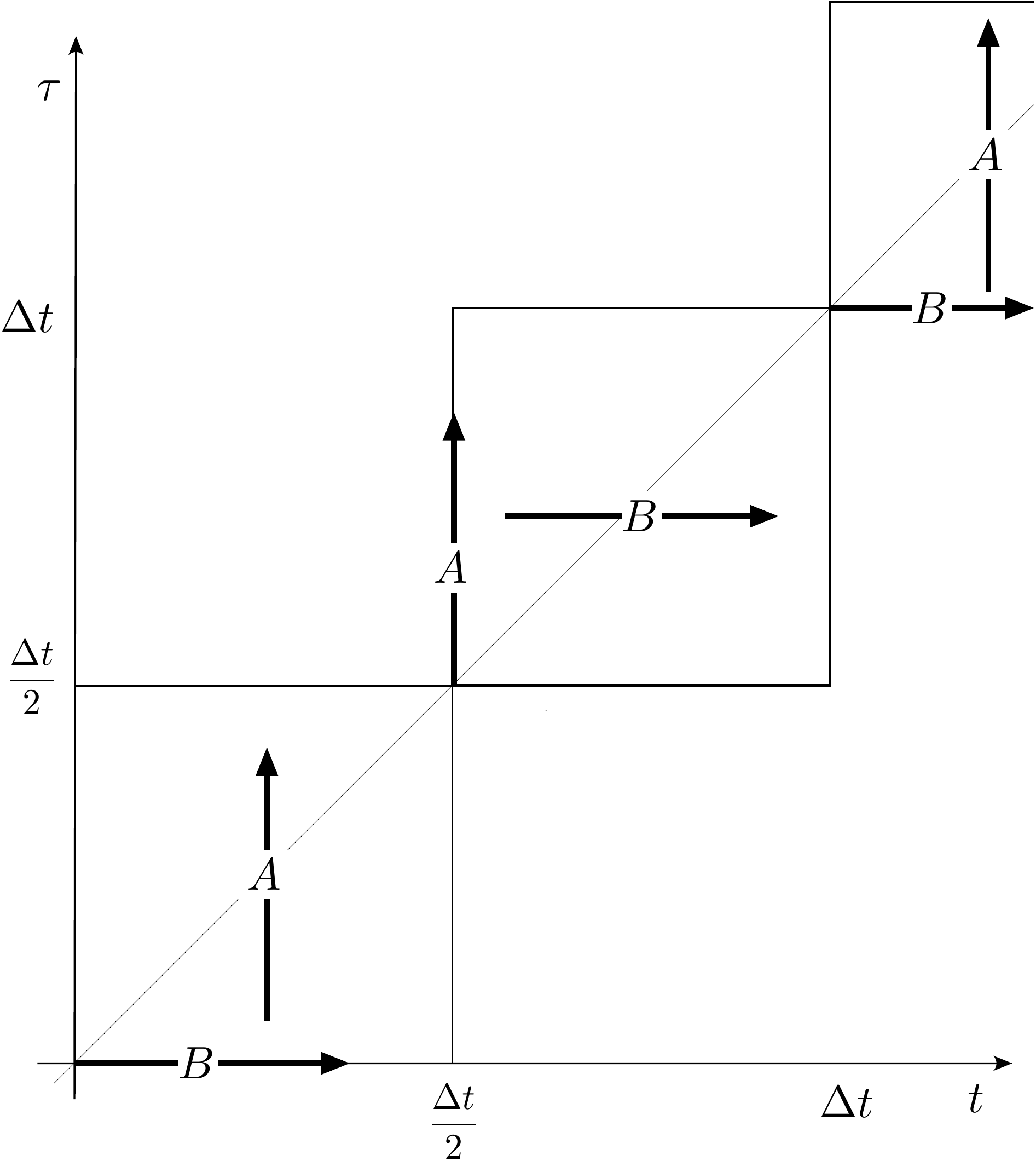}
  \caption{A schematic view of Strang splitting.}
  \label{fig:StrangSplit}
\end{figure}

The aim is now to show that
\begin{equation*}
        \norm{v\left(t,t\right)-u\left(t\right)}=\Oh\left(\Dt^2\right)
\end{equation*}
in an appropriate norm.  Here the abstract analysis 
of Section \ref{subsec:tdouble} applies, and we find:
\begin{description}
\item[(i)] On the domains $[t_{n},t_{n+1/2}]\times[t_{n},t_{n+1/2}]$ we
  have
  \begin{align}
      w_t-dC(u)[w]&=F(t)+ \int_0^1 (1-\alpha) d^{(2)}C(u+\alpha w)[w]^2 d\alpha,\notag\\
      F(t,\tau)&=v_t(t,\tau)-B(v(t,\tau)).\label{eq:w_eq2AA}
    \end{align}
 The forcing term satisfies the following time development
  \begin{equation*}
    \begin{aligned}
      F_\tau&=dA(v)[F]+[A,B](v),\\
      F_t&=v_{tt}-dB(v)[F]-dB(v)[B(v)].
    \end{aligned}
  \end{equation*}

\item[(ii)] On the domains $[t_{n+1/2},t_{n+1}]\times[t_{n+1/2},t_{n+1}]$
  we have (here we write $w_t$ rather than $w_\tau$ since $w$ is a
  function of one variable only)
  \begin{align}
      w_t- dC(u)[w]&=G(t)+ \int_0^1 (1-\alpha) d^{(2)}C(u+\alpha
      w)[w]^2 d\alpha, \notag\\ 
      G(t,\tau)&=v_\tau(t,\tau)-A(v(t,\tau)). \label{eq:w_eq2aA}
    \end{align}
  The forcing term satisfies the following time development
  \begin{equation*}
    \begin{aligned}
      G_t&=dB(v)[G]+[B,A](v),\\
      G_\tau&=v_{\tau\tau}-dA(v)[G]-dA(v)[A(v)].
    \end{aligned}
  \end{equation*}
\end{description}
We extend $F$ and $G$ to all of $ \widetilde\Omega_{\Dt}$ using the same definitions, \eqref{eq:w_eq2AA} and \eqref{eq:w_eq2aA}, respectively. Observe that this implies that
$F=0$ on  $[t_{n+1/2},t_{n+1}]\times[t_{n+1/2},t_{n+1}]$, while $G=0$ on $[t_{n},t_{n+1/2}]\times[t_{n},t_{n+1/2}]$.
The total forcing term is defined by
\begin{equation*}
  H(t,\tau)=F(t,\tau)+G(t,\tau).
\end{equation*}

\subsection{Ordinary differential equations}
\label{subsec:ode_strang2}
One can consider the case of ordinary differential equations, as we did in 
Subsection \ref{subsec:odes} for Godunov splitting, but for reasons of 
brevity we will only revisit the example in  Remark \ref{rem:ODE1}.
\begin{remark}\label{rem:ODE2}

We find, using the definitions
  \eqref{eq:GodunovDef}, that 
\begin{equation}
  \label{eq:ex1A_ODE}
  v(t,\tau)=\begin{cases}
\frac{v(t_n,t_n)e^{-(\tau-t_n)}}{1-v(t_n,t_n)(t-t_n)}, & 
\text{for $t,\tau\in[t_n,t_{n+1/2}]$}, \\
\frac{v(t_{n+1/2},t_{n+1/2})}{e^{\tau-t_{n+1/2}}-
v(t_{n+1/2},t_{n+1/2})(t-t_{n+1/2})}, & 
\text{for $t,\tau\in[t_{n+1/2},t_{n+1}]$},
\end{cases}
\end{equation}
for $n=0,\dots,\lfloor T/\Dt \rfloor$. By induction we determine
\begin{equation}
  \label{eq:ex1Aa_ODE}
  v(t_n,t_n)= \frac{u_0(1-e^{-\Dt})}{(1-e^{-\Dt})e^{t_n}
+u_0\Dt(e^{t_n}-1)(e^{\Dt}+1)/2}.
\end{equation}
In this case we compute,
when we for convenience write  $\alpha=u_0(1-e^{t_n})$, 
\begin{equation}
  \label{eq:exAA_ODE}
  \begin{aligned}
    \abs{w(t_n)}&=\abs{v(t_n,t_n)-u(t_n)}\\
    &=\abs{\frac{u_0(1-e^{-\Dt})}{(1-e^{-\Dt})e^{t_n}
        +u_0\Dt(e^{t_n}-1)(e^{\Dt}+1)/2}
      -\frac{u_0 }{e^{t_n}+u_0(1-e^{t_n})}}\\
    &= \abs{u_0} \abs{\frac{1}{e^{t_n}+\alpha
        \Dt(1+e^{\Dt})/(2(1-e^{-\Dt}))}
      -\frac{1}{e^{t_n}+\alpha}}\\
    &\le \frac{\abs{u_0}^2 e^{t_n}(e^{t_n}-1)}{\abs{(e^{t_n}+\alpha
        \Dt(1+e^{\Dt})/(2(1-e^{-\Dt})))(e^{t_n}+\alpha)}}\,
    \abs{1- \frac{\Dt(1+e^{\Dt})}{2(1-e^{-\Dt})}}\\
    &\le \Oh(\Dt^2).
  \end{aligned}
\end{equation}
\begin{figure}[tbp]
  \centering \includegraphics[width=0.45\linewidth]{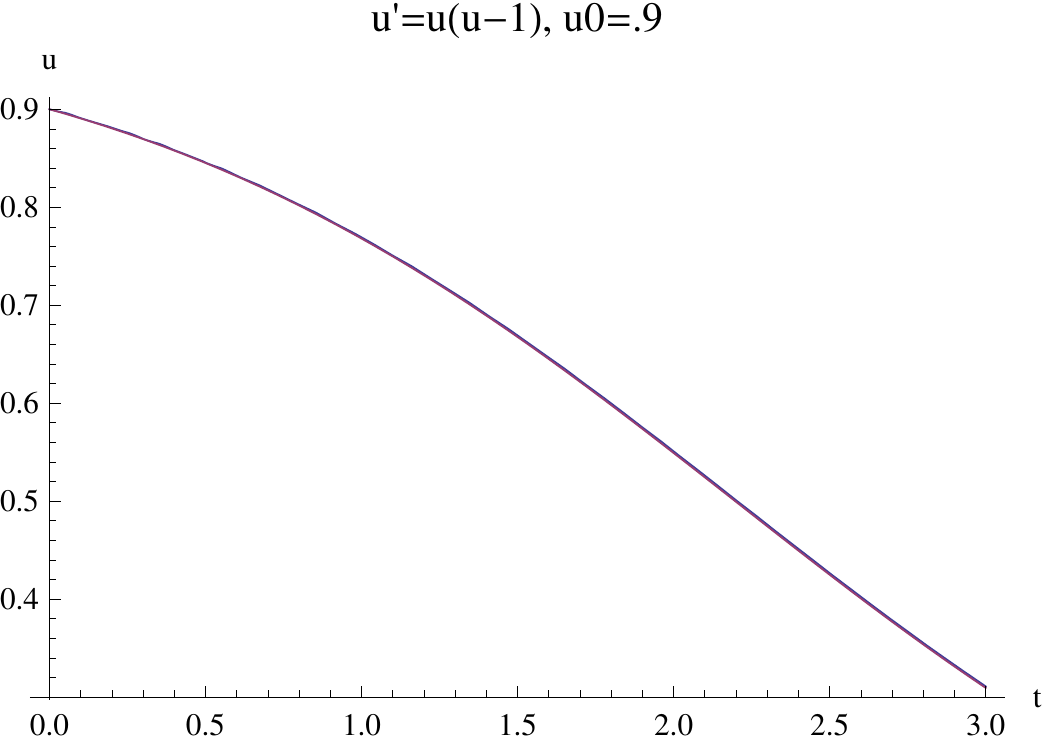}
  \includegraphics[width=0.45\linewidth]{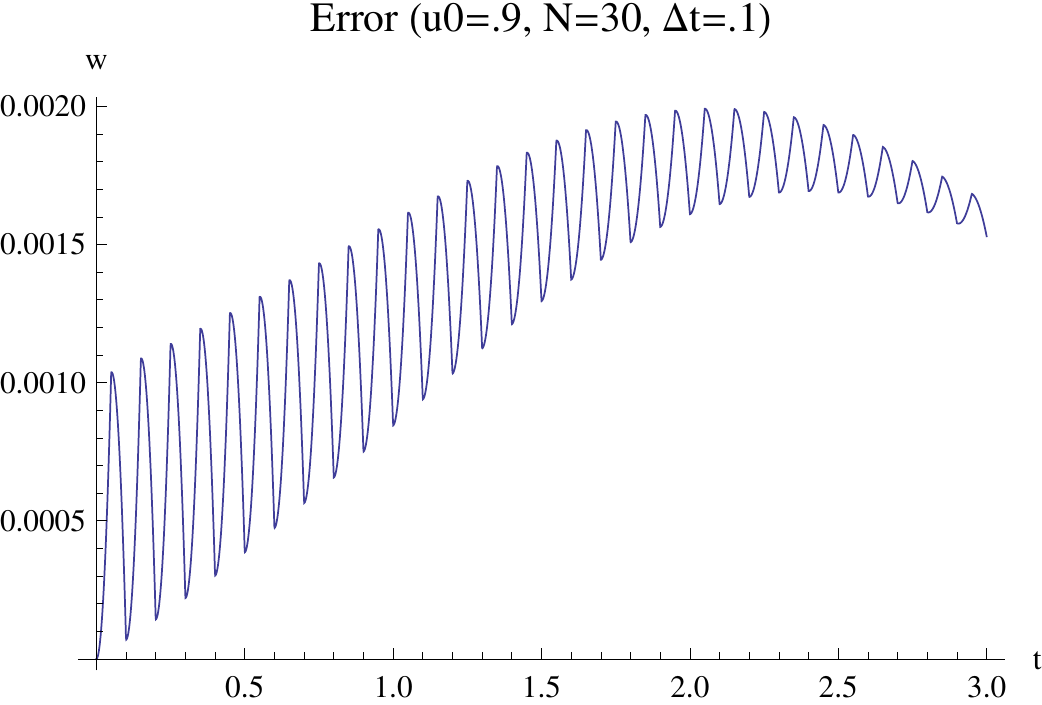}
      \caption{(Left) The exact  (blue) and the approximate (red)
        solution. (Right) The error $v(t,t)-u(t)$.}
  \label{fig:odeSplit1}
\end{figure}
The example is illustrated in Figure \ref{fig:odeSplit1}. Observe the
strong oscillations in the error, it is these oscillations which
prevent the error from growing too large.
\end{remark}

\subsection{The KdV equation}
\label{subsec:kdv_strang2}
For the KdV equation, we use $B(v)=vv_x$ and $A(v)=-v_{xxx}$, and the
above analysis yields:
\begin{description}
\item[(i)]
  On domains $[t_{n},t_{n+1/2}]\times[t_{n},t_{n+1/2}]$ we have
  \begin{equation}
    \label{eq:w_eq2A}
    \begin{aligned}
      w_t-(uw)_x+w_{xxx}&=F(t)+ww_x,\\
      F(t,\tau)&=v_t-vv_x.
    \end{aligned}
  \end{equation}
  The forcing term satisfies the following time development
  \begin{equation}
    \label{eq:F_tau2A}
    \begin{aligned}
      F_\tau&=-F_{xxx}-\frac32\partial_x^2(v_x)^2,\\
      F_t&=v_{tt}-(vF)_x-(2vv_x^2+v^2v_{xx}),
    \end{aligned}
  \end{equation}
  since $dB(v)[B(v)]=v(vv_x)_x+v_x(vv_x)=2vv_x^2+v^2v_{xx}$.\\

\item[(ii)]
  On domains $[t_{n+1/2},t_{n+1}]\times[t_{n+1/2},t_{n+1}]$ we have
  \begin{equation}
    \label{eq:w_eq2B}
    \begin{aligned}
      w_t-(uw)_x+w_{xxx} & =G(t)+ww_x,\\
      G(t,\tau)&=v_\tau(t,\tau)+v_{xxx}.
    \end{aligned}
  \end{equation}
  The forcing term satisfies the following time development
  \begin{equation}
    \label{eq:F_tau2B}
    \begin{aligned}
      G_t&=(vG)_x+\frac32\partial_x^2(v_x)^2,\\
      G_\tau&=v_{\tau\tau}+G_{xxx}-v^{(6)},
    \end{aligned}
  \end{equation}
  since $dA(v)[A(v)]=(v_{xxx})_{xxx}$.
\end{description}

To start the bootstrapping procedure we fix an odd integer $s$ and a
positive constant $\alpha$, whose values will be determined in the course
of the argument. Now assume that 
$$
\norm{v(t,\tau)}_{H^{\hat{k}}} \le \alpha, \quad (t,\tau)\in\widetilde\Omega_\Dt,
$$
where $\hat{k}=(s-1)/2$. As a consequence, we have
$$
\norm{v(t,\tau)}_{H^s}\le K_\alpha, \quad (t,\tau)\in\widetilde\Omega_\Dt,
$$
cf.~the proof of Lemma \ref{lem:indu}, which 
can be easily adapted to Strang splitting.

We need to introduce the function
\begin{equation}  \label{eq:z1}
z(t)=w(t)+w(t+\frac{\Dt}{2})=w(t)+\tilde w(t).
\end{equation}
In the following we will write $\tilde\phi(t)=\phi(t+\frac{\Dt}{2})$
for any function $\phi$.  Straightforward calculations yield that $z$
satisfies
\begin{equation}  \label{eq:z2}
z_t-(\frac12 z^2+uz-z_{xx})_x=H+\widetilde H+\big(\tilde w(\tilde u-u)-w\tilde w \big)_x.
\end{equation}

By our techniques, we must work in $H^{s-9}$ (so at least $s\ge 9$), therefore set 
\begin{equation*}
E(t)=\norm{z(t)}_{H^{s-9}(\R)}.
\end{equation*}
By taking the $H^{s-9}$ inner product with $z$ in \eqref{eq:z2}
 we get
\begin{align*}
    \frac{1}{2}\frac{d}{dt} E^2(t)&=(z,(\frac12 z^2+uz-z_{xx})_x)_{H^{s-9}}+ (z,H+\widetilde H)_{H^{s-9}}
    +(z,(\tilde w(\tilde u-u)-w\tilde w)_x)_{H^{s-9}} \\
&=(z,zz_x+u_xz+uz_x)_{H^{s-9}}+ (z,H+\widetilde H)_{H^{s-9}}
    -(z_x,\tilde w(\tilde u-u)-w\tilde w)_{H^{s-9}} \\
&\le  (z,zz_x+u_xz+uz_x)_{H^{s-9}}  +E(t) \norm{H+\widetilde H}_{H^{s-9}}
+E(t) \norm{\tilde w(\tilde u-u)-w\tilde w}_{H^{s-9}}
    \\    
    &\le K_\alpha E^2(t)+K E(t)\big( \norm{H+\widetilde H}_{H^{s-9}}+\norm{\tilde w(\tilde u-u)-w\tilde w}_{H^{s-9}} \big).
\end{align*}
Thus by Gronwall's inequality
\begin{equation}
\begin{aligned}
  \label{eq:EStrang}
  E(t)&\le E(0)+\int_0^t e^{K_\alpha (t -\sigma)}
  \big( \norm{H(\sigma,\sigma)+\widetilde H(\sigma,\sigma)}_{H^{s-9}}\\
  &\qquad\qquad\qquad +\norm{\tilde w(\sigma)(\tilde
  u(\sigma)-u(\sigma))-w(\sigma)\tilde w(\sigma)}_{H^{s-9}}
  \big)\,d\sigma\\ 
  &\le E(0)+ e^{K_\alpha t} \int_0^t
  \big( \norm{H(\sigma,\sigma)+\widetilde H(\sigma,\sigma)}_{H^{s-9}}\\
  &\qquad\qquad\qquad\qquad +\norm{\tilde w(\sigma)(\tilde
    u(\sigma)-u(\sigma))-w(\sigma)\tilde w(\sigma)}_{H^{s-9}}
  \big)\,d\sigma.
  \end{aligned}
\end{equation}

Next we turn to the detailed estimate of each of the terms in
\eqref{eq:EStrang}. We start with most involved one, the forcing term,
which can be estimated as follows. We consider the term $F$ first.
Since $F(t,t_n)=0$, we easily see that $F_t(t_n,t_n)=0$; thus (cf.~
\eqref{eq:F_tau2A})
$(F_t+F_\tau)(t_n,t_n)=-\frac32\partial_x^2(v_x)^2(t_n,t_n)$. Thus
\begin{align*}
  F(t,t)&=(v_t-vv_x)(t,t) \\
  &=F(t_n,t_n)+(F_t+F_\tau)(t_n,t_n)(t-t_n) \\
  &\quad +\frac{(t-t_n)^2}2\int_0^1(F_{tt}+2F_{t\tau}+F_{\tau\tau})
  (\sigma(t- t_n)+t_n,\sigma(t- t_n)+t_n)d\sigma
  \\
  &= -\frac32\partial_x^2(v_x)^2(t_n,t_n)(t-t_n) \\
  &\quad +\frac{(t-t_n)^2}2\int_0^1(F_{tt}+2F_{t\tau}+F_{\tau\tau})
  (\sigma(t- t_n)+t_n,\sigma(t- t_n)+t_n)d\sigma.
\end{align*}
 As for the second derivatives, we find
\begin{equation}
  \label{eq:force5b}
  \begin{aligned}
    F_{\tau\tau}&=-F_{\tau xxx}-\frac32\partial^2_x\partial_\tau
    (v_x)^2\\
&=\partial^6_x F+\frac32\partial^5_x(v_x)^2-3\partial^2_x(v_x
v^{(6)}),\\[2mm]
 F_{\tau t}&=-v_{ttxxx}-(v_{xxx}v_x+vv_{xxxx})_t, \\[2mm]
 F_{tt}&=v_{ttt}-v_{tt}v_x-2v_tv_{xt}-vv_{xtt}.
  \end{aligned}
\end{equation}
Similarly we find for the forcing term $G$ the following
estimates. Since $G(t_{n+1/2},\tau)=0$, we easily see that
$G_\tau(t_{n+1/2},t_{n+1/2})=0$; thus (cf.~\eqref{eq:F_tau2B})
$(G_t+G_\tau)(t_{n+1/2},t_{n+1/2})=\frac32\partial_x^2(v_x)^2(t_{n+1/2},t_{n+1/2})$. Thus 
\begin{align*}
  G(t+\frac\Dt2,t+\frac\Dt2)&=(v_\tau+v_{xxx})(t+\frac\Dt2,t+\frac\Dt2) \\
  &=G(t_{n+1/2},t_{n+1/2})+(G_t+G_\tau)(t_{n+1/2},t_{n+1/2})(t-t_n)\\
  &\quad +\frac{(t-t_n)^2}2\int_0^1(G_{tt}+2G_{t\tau}+G_{\tau\tau})
  (\sigma(t- t_n)+t_n,\sigma(t- t_n)+t_n)d\sigma \\
  &= \frac32\partial_x^2(v_x)^2(t_{n+1/2},t_{n+1/2})(t-t_n) \\
  &\quad +\frac{(t-t_n)^2}2\int_0^1(G_{tt}+2G_{t\tau}+G_{\tau\tau})
  (\sigma(t- t_n)+t_n,\sigma(t- t_n)+t_n)d\sigma.
\end{align*}
The second derivatives $G_{tt}+2G_{t\tau}+G_{\tau\tau}$ will have to
be considered similarly to those for $F$. 
These read 
\begin{align*}
  G_{tt}&=2v_x^2 G + 4 v v_x G_x + v^2G_{xx} + G v v_{xx} \\
  &\qquad + 6 v_{xx}\left(v_x v\right)_{xx} + v v_x v_{xxx} +
  v\left(v v_x\right)_{xxx},\\
  G_{\tau t} &= v_{\tau\tau} v_x + 2v_\tau v_{\tau x} + v v_{\tau\tau
    x} + \partial^4\left(Gv\right) + \frac{3}{2} \partial^5\left(v_x^2\right)
  - \partial^6\left(vv_x\right), \\
  G_{\tau\tau} &= v_{\tau\tau\tau} + v_{\tau\tau x x x}.
\end{align*}
\begin{lemma}
  \label{lem:FTlemma} We have the estimate 
  \begin{equation*}
    \norm{\Phi(t,\tau)}_{H^{s-9}}^2+ \norm{\Psi(\bar t,\bar
      \tau)}_{H^{s-9}}^2\le K_\alpha,  
    \quad (t,\tau),(\bar t,\bar\tau) \in\widetilde\Omega_\Dt,
  \end{equation*}
  where $\Phi=(F_{tt}+2F_{t\tau}+F_{\tau\tau})$ and 
  $\Psi=(G_{tt}+2G_{t\tau}+G_{\tau\tau})$.
\end{lemma}
\begin{proof}
We have that $\norm{F}_{H^{s-9}} \le K_\alpha\Dt$, we shall get a similar
estimate for $\norm{G}_{H^{s-9}}$. For $t_n\le \tau \le t_{n+1/2}$,
we have that 
$$
G_t = v_x G + v G_x + \frac{3}{2}\partial^2 \left(v_x^2\right),
$$
and $G(t,t_n)=0$. Taking the $H^{s-9}$ inner product with $G$ we get
\begin{equation*}
  \frac{1}{2}\frac{d}{dt}\norm{G}_{H^{s-9}}^2 = 
  \sum_{k=0}^{s-9} \int \partial^k\left(v_x G\right)\partial^k G
  + \partial^k \left(v G_x\right)\partial^k G +
  \frac{3}{2}\partial^{k+2}\left( v_x^2\right)\partial^k G \,dx.
\end{equation*}

The first term expands by the Leibniz rule; a typical term in this
expansion reads (here $0\le j\le k\le s-9$),
$$
\int \partial^{j+1} v \partial^{k-j} G \partial^k G \,dx \le
\norm{\partial^{j+1} v}_{L^\infty} \norm{G}_{H^{s-9}}^2.
$$

Similarly the second term can be expanded and estimated, except the
term containing $\partial^{k+1}G$, which is estimated as
\begin{align*}
        \int v \partial^{k+1}G\partial^k G \,dx 
         = \frac{1}{2} \int v \partial\left(\partial^k G\right)^2 \,dx 
        & = -\frac{1}{2} \int\partial v \left(\partial^k G\right)^2 \,dx
        \\ & \le \frac12 \norm{\partial v}_{L^\infty} \norm{G}_{H^{s-9}}^2.
\end{align*}

We bound the last term as
$$
\Bigl| \int \partial^{k+2}\left(v_x^2\right) \partial^k G \,dx\Bigr|
\le 
K \norm{v}_{H^{s}}^2 \norm{G}_{H^{s-9}}.
$$  

Summing up, we get
$$
\frac{d}{dt}\norm{G}_{H^{s-9}}^2 \le 
K_\alpha \norm{G}_{H^{s-9}}^2 + K_\alpha \norm{G}_{H^{s-9}}.
$$
Using Gronwall's inequality and that $G(t_n,\tau)=0$ we get
\begin{equation*}
        \norm{G}_{H^{s-9}}\le K_\alpha \Dt.
\end{equation*}

We also need estimates for $v_t$, $v_{tt}$ and $v_{ttt}$ where
$t,\tau\in [t_n,t_{n+1/2}]$. In this set,
$v_\tau=-v_{xxx}$, and this evolution 
preserves the $H^k$ norm. For $\tau=t_n$, 
\begin{align*}
  v_t&=vv_x, \\
  v_{tt} &= v_t v_x + v v_{tx},\\
  v_{ttt} &= v_{tt} v_x + 2 v_t v_{tx} + v v_{ttx}.
\end{align*}
Therefore
\begin{equation*}
        \begin{aligned}
                \norm{v_t}_{H^k} &\le \norm{v}_{H^k}\norm{v}_{H^{k+1}}
                \le K\norm{v}_{H^{k+1}}^2,\\ 
                \norm{v_{tt}}_{H^{k}} &\le K \norm{v}_{H^{k+2}}^3,\\
                \norm{v_{ttt}}_{H^k} &\le K \norm{v}_{H^{k+3}}^4.
        \end{aligned}
\end{equation*}
Next we turn to estimates of $v_\tau$, $v_{\tau\tau}$ and 
$v_{\tau\tau\tau}$ in the set $t,\tau\in [t_{n+1/2},t_{n+1}]$. 
Here $v_t=vv_x$, and setting $\omega=v_\tau$, $\theta=v_{\tau\tau}$ 
and $\gamma=v_{\tau\tau\tau}$, we get
\begin{align*}
  \omega_t &=  v_x \omega + v \omega_x ,
  &\omega(t_{n+1/2},\tau)=-\partial^3 v,\\
  \theta_t &= 2 \omega \omega_x + v_x\theta + v \theta_x, 
  &\theta(t_{n+1/2},\tau)=\partial^6 v,\\
  \gamma_t &= 3\theta \omega_x + 3 \theta_x \omega + v_x \gamma+v\gamma_x,  
  &\gamma(t_{n+1/2},\tau)=-\partial^9 v.
\end{align*}
Starting with $\omega$, for $k<s-9$ we get 
$$
\frac{1}{2}\frac{d}{dt}\norm{\omega}_{H^k}^2 = \sum_{j=0}^k
\int \partial^j \left(v_x \omega\right) \partial^j \omega + \partial^j
\left(v\omega_x\right) \partial^j \omega \,dx.
$$
Using Leibniz' rule, all terms except the ultimate one will be of the type 
$$
\Bigl| \int \partial^\ell v \partial^{j-\ell} \omega \partial^j \, dx\Bigr| 
\le \norm{v}_{H^s} \norm{\omega}_{H^k}^2.
$$
The ultimate term (with one too many derivatives on $\omega$) is as usual
estimated as
$$
\Bigl| \int v \partial^{k+1}\omega \partial^k \omega \,dx\Bigr| = 
\frac{1}{2}\Bigl| \int \partial v\big(\partial^k \omega\big)^2\,dx\Bigr| 
\le \norm{v}_{H^s} \norm{\omega}_{H^k}^2.
$$
Gronwall's inequality then yields
$$
\norm{\omega}_{H^k}(t)\le e^{K_\alpha \Dt} \norm{\omega}_{H^k}(t_{n+1/2})\le 
K_\alpha \norm{v}_{H^{k+3}}\left(t_{n+1/2},\tau\right).
$$
Reasoning similarly for $\theta$, we find that
\begin{align*}
\norm{\theta}_{H^k}(t)&\le e^{K_\alpha \Dt}\left(
\max_{s\in [t_{n+1/2},t_{n+1}]} \norm{\omega}_{H^{k+1}}^2(s) +
\norm{\theta}_{H^k}(t_{n+1/2}) 
\right)\\
&\le K_\alpha \left(\norm{v}_{H^{k+4}}^2 +
  \norm{v}_{H^{k+6}}\right)\left(t_{n+1/2},\tau\right).
\end{align*}
Finally, the estimate for $\gamma$ reads
\begin{align*}
  \norm{\gamma}_{H^k}&\le
  K_\alpha\left(\norm{\theta}_{H^{k+1}}\norm{\omega}_{H^{k+1}} +
    \norm{\gamma}_{H^k} \right)\left(t_{n+1/2},\tau\right)\\
  &\le K_\alpha\left(\norm{v}_{H^{k+4}}^2 + \norm{v}_{H^{k+5}}^4 +
      \norm{v}_{H^{k+7}}^2 + \norm{v}_{H^{k+9}}\right) \left(t_{n+1/2},\tau\right).
\end{align*}
Summing up
\begin{equation*}
  \norm{v_\tau}_{H^k} \le K_\alpha\norm{v}_{H^{k+3}},\,
  \norm{v_{\tau\tau}}_{H_k}\le K_\alpha \norm{v}_{H^{k+6}},\,
  \norm{v_{\tau\tau\tau}}_{H_k}\le K_\alpha \norm{v}_{H^{k+9}}. 
\end{equation*}
Now
$$
\norm{\Phi}_{H^{s-9}} \le \norm{F_{\tau\tau}}_{H^{s-9}} +
2\norm{F_{t\tau}}_{H^{s-9}} + \norm{F_{tt}}_{H^{s-9}}.
$$

Working in the square $[t_n,t_{n+1/2}]^2$, 
from \eqref{eq:force5b}, 
\begin{align}
  \norm{F_{\tau\tau}}_{H^{s-9}} &\le 
  \norm{\partial^6 F}_{H^{s-9}} + \frac{3}{2}
  \norm{\partial^5\left(\partial v\right)^2}_{H^{s-9}} + 
  3\norm{\partial^2\left(\partial v\partial^6v\right)}_{H^{s-9}}\notag\\
  &\le \norm{F}_{H^{s-3}} + K_\alpha \norm{v}_{H^{s-3}}^2 + K_\alpha
  \norm{v}_{H^{s-1}}^2\label{eq:Ftautau}\\
  &\le K_\alpha,\notag
\end{align}
and 
\begin{align}
  \norm{F_{t\tau}}_{H^{s-9}} &\le \norm{v_{tt}}_{H^{s-6}} + K_\alpha\Bigl(
  \norm{v}_{H^{s-8}}\norm{v_t}_{H^{s-6}} +
  \norm{v_t}_{H^{s-9}}\norm{v}_{H^{s-6}} \notag 
  \\&\hphantom{\le \norm{v_{tt}}_{H^{s-6}} + K_\alpha\Bigl(}\quad +
  \norm{v}_{H^{s-9}} \norm{v_t}_{H^{s-5}} +
  \norm{v_t}_{H^{s-9}}\norm{v}_{H^{s-5}}\Bigr)\notag\\
  &\le K_\alpha\Bigl(
  \norm{v}_{H^{s-4}} + \norm{v}_{H^{s-8}}\norm{v}_{H^{s-5}}  
   \label{eq:Fttau}\\
  &\hphantom{\le K_\alpha \Bigl(}\quad+ 
  \norm{v}_{H^{s-7}} \norm{v}_{H^{s-7}} +
  \norm{v}_{H^{s-9}}\norm{v}_{H^{s-4}} \notag \\
  &\hphantom{\le K_\alpha \Bigl(}\quad +
  \norm{v}_{H^{s-8}}\norm{v}_{H^{s-5}}\Bigr)\notag\\
  &\le K_\alpha, \notag
\end{align}
and also
\begin{align}
  \norm{F_{tt}}_{H^{s-9}} & \le \norm{v_{ttt}}_{H^{s-9}} +
  \norm{v_{tt}v_x}_{H^{s-9}}
  +2\norm{v_tv_{xx}}_{H^{s-9}} + \norm{v v_{txx}}_{H^{s-9}}\notag\\
  &\le K_\alpha\Bigl(
  \norm{v}_{H^{s-6}} + \norm{v}_{H^{s-7}}\norm{v}_{H^{s-8}} +
  \norm{v}_{H^{s-9}}\norm{v}_{H^{s-6}}\Bigr) \notag\\
  &\le K_\alpha. \label{eq:Ftt}
\end{align}
Hence $\norm{\Phi}_{H^{s-9}}\le K_\alpha$. 

In the second square $[t_{n+1/2},t_n]^2$, we write 
$$
\norm{\Psi}_{H^{s-9}} \le \norm{G_{tt}}_{H^{s-9}} +
2\norm{G_{t\tau}}_{H^{s-9}} + \norm{G_{\tau\tau}}_{H^{s-9}}.
$$
Each term above is estimated individually as
\begin{align}
  \norm{G_{tt}}_{H^{s-9}}&\le 
  2\norm{v_x^2 G}_{H^{s-9}} + 4\norm{v v_x G_x}_{H^{s-9}} +
  \norm{v^2 G_{xx}}_{H^{s-9}} \notag\\
  &\qquad 
  +3\norm{v_{xx} \left(v^2\right)_{xx}}_{H^{s-9}} + \norm{v v_{xx}
    v_{xxx}}_{H^{s-9}} +
  \norm{v\left(vv_x\right)_{xxx}}_{H^{s-9}}\notag\\
  &\le K_\alpha\Bigl( \norm{v}_{H^{s-8}} \norm{G}_{H^{s-9}} +
  \norm{v}_{H^{s-8}}\norm{G}_{H^{s-8}} + \norm{v}_{H^{s-9}}
  \norm{G}_{H^{s-7}}\notag\\
  &\hphantom{K_\alpha\Bigl(}\quad + 
  \norm{v}_{H^{s-7}}^2 + \norm{v}_{H^{s-6}}^3 +
  \norm{v}_{H^{s-5}}^3\Bigr)\notag\\
  &\le K_\alpha, \label{eq:Gtt} \\[2mm]
  \norm{G_{t\tau}}_{H^{s-9}}&\le \norm{v_{\tau\tau} v_x}_{H^{s-9}} +
  2\norm{v_\tau v_{xx}}_{H^{s-9}} + \norm{v v_{\tau\tau x}}_{H^{s-9}}
  + \norm{\partial^4(Gv)}_{H^{s-9}}\notag\\
  &\qquad +
  \frac{3}{2}\norm{\partial^5(v_x^2)}_{H^{s-9}}
  +\norm{\partial^6(vv_x)}_{H^{s-9}}\notag\\
  &\le K_\alpha\Bigl(\norm{v}_{H^{s-3}} \norm{v}_{H^{s-8}} +
  \norm{v}_{H^{s-6}} \norm{v}_{H^{s-8}} +
  \norm{v}_{H^{s-9}} \norm{v}_{H^{s-1}} \notag\\
  &\hphantom{\le K_\alpha\Bigl(}\quad +
  \norm{G}_{H^{s-5}} \norm{v}_{H^{s-5}} + 
  \norm{v}_{H^{s-3}}^2 + 
  \norm{v}_{H^{s-3}}  \norm{v}_{H^{s-2}}\Bigr)  \notag\\
  &\le K_\alpha, \label{eq:Gttau} \\[2mm]
  \norm{G_{\tau\tau}}_{H^{s-9}} &\le \norm{v_{\tau\tau\tau}}_{H^{s-9}}
  + \norm{v_{\tau\tau x x x}}_{H^{s-9}} \notag\\
  &\le K_\alpha \Bigl( \norm{v}_{H^s} + \norm{v}_{H^s}\Bigr)\notag\\
  &\le K_\alpha. \label{eq:Gtautau}
\end{align}
Collecting \eqref{eq:Ftt}, \eqref{eq:Fttau}, \eqref{eq:Ftautau},
\eqref{eq:Gtt}, \eqref{eq:Gttau} and \eqref{eq:Gtautau} finishes the
proof of the lemma.
\end{proof}
\begin{lemma}\label{lem:vxlip} The map
$$
[0,T]\ni t\mapsto \partial^2 \big(v_x^2\big)(t,t) \in H^{s-9},
$$
is Lipschitz continuous with Lipschitz constant  at most $K_\alpha$.
\end{lemma}
\begin{proof}
Set $w(t,\tau)=\partial^2(v_x^2)(t,\tau)$. Then we have
\begin{align*}
        w_t &= 4 v_{xx} v_{txx} + 2 v_{txxx}v_x + 2 v_{xxx}v_{tx},\\
        w_\tau &= 4 v_{xx} v_{\tau xx} + 2 v_{\tau xxx}v_x + 2
        v_{xxx}v_{\tau x}.
\end{align*}
In the square $[t_n,t_{n+1/2}]^2$ we have $v_\tau=-v_{xxx}$, thus
$$
w_\tau=-4 v_{xx} \partial^5 v - 2 \partial^6 v v_x - 2
v_{xxx}\partial^4v.
$$
Hence in this square
\begin{align*}
        \norm{w_t+w_\tau}_{H^{s-9}} &\le
        K\Bigl( \norm{v}_{H^{s-7}} \norm{v_t}_{H^{s-7}} +
        \norm{v_t}_{H^{s-6}}\norm{v}_{H^{s-8}} \\
        &\hphantom{\le C\Bigl)}\quad+
        \norm{v}_{H^{s-6}}\norm{v_t}_{H^{s-8}}
        + \norm{v}_{H^{s-7}} \norm{v}_{H^{s-4}}
        \\
        &\hphantom{\le C\Bigl)}\quad+
        \norm{v}_{H^{s-3}} \norm{v}_{H^{s-1}} +
        \norm{v}_{H^{s-6}}\norm{v}_{H^{s-5}} \Bigr)\\
        &\le K_\alpha.
\end{align*}
In the second square $[t_{n+1/2},t_{n+1}]^2$ we have $w_t=vv_x$, and
$$
w_t = 4 v_{xx}\left(vv_x\right)_{xx} +
2\left(vv_x\right)_{xxx} v_x +
2v_{xxx} \left(vv_x\right)_{x}.
$$
Therefore in this square
\begin{align*}
        \norm{w_t+w_\tau}_{H^{s-9}} &\le
        K\Bigl( \norm{v}_{H^{s-7}} \norm{vv_x}_{H^{s-7}} +
        \norm{vv_x}_{H^{s-6}} \norm{v}_{H^{s-9}} \\
        &\hphantom{\le C\Bigl(}\quad+
        \norm{v}_{H^{s-6}} \norm{vv_x}_{H^{s-8}} +
        \norm{v}_{H^{s-8}}\norm{v_\tau}_{H^{s-8}}
        \\
        &\hphantom{\le C\Bigl(}\quad+
        \norm{v}_{H^{s-8}}\norm{v_\tau}_{H^{s-6}} +
        \norm{v}_{H^{s-6}}\norm{v_\tau}_{H^{s-8}}\Bigr)\\
        &\le
        K\Bigl( \norm{v}_{H^{s-7}}^2 \norm{v}_{H^{s-6}} +
        \norm{v}_{H^{s-7}}\norm{v}_{H^{s-6}} \norm{v}_{H^{s-9}}
        \\
        &\hphantom{\le C\Bigl(}\quad+
        \norm{v}_{H^{s-6}} \norm{v}_{H^{s-8}}\norm{v}_{H^{s-7}} +
        \norm{v}_{H^{s-8}}\norm{v}_{H^{s-5}}
        \\
        &\hphantom{\le C\Bigl(}\quad+
        \norm{v}_{H^{s-8}}\norm{v}_{H^{s-3}} +
        \norm{v}_{H^{s-6}}\norm{v}_{H^{s-5}}\Bigr)\\
        &\le K_\alpha.
\end{align*}
\end{proof}
From this lemma it follows that 
\begin{equation*}
        \Gamma(t)=-\partial^2(v_x)^2(t,t)
        +\partial^2(v_x)^2\left(t+\frac{\Dt}{2},t+\frac{\Dt}{2}\right) 
\end{equation*}
satisfies
\begin{equation*}
        \norm{\Gamma(t)}_{H^{s-9}}\le K_\alpha \Dt,
        \quad t\in[0,T].
\end{equation*}

The following lemma will be convenient:
\begin{lemma}\label{lem:cancellation}
For $t\in [t_m,t_{m+1/2}]$ we have
\begin{align}
  \sum_{j=0}^{s-9}
  \int\left(\partial^j\left(F(t)+G\left(t+\Dt/2\right)\right)\right)^2\,dx
  &\le K_\alpha \Dt^4,\label{eq:FGone}\\
  \sum_{j=0}^{s-9}
  \int\left(\partial^j\left(F(t)+G\left(t-\Dt/2\right)\right)\right)^2\,dx
  &\le K_\alpha \Dt^4.\label{eq:FGtwo}
\end{align}
\end{lemma}

\begin{proof}
We show \eqref{eq:FGone}, \eqref{eq:FGtwo} is proved similarly. 
By a Taylor expansion, for $t\in [t_m,t_{m+1/2}]$
\begin{align*}
  \sum_{j=0}^{s-9}&
  \int\left(\partial^j\left(F(t)+G\left(t+\Dt/2\right)\right)\right)^2\,dx
  \\
   &=\sum_{j=0}^{s-9} \int
   \Bigl[\partial^j \Bigl(
   -\frac32\partial^2\left(v_x\right)^2\left(t_m\right)\left(t-t_m\right)
   +\frac{(t-t_m)^2}{2}\int_0^1\Phi\left(\sigma(t-t_m)+t_m\right)d\sigma
   \\
   &\hphantom{=\sum_{j=0}^{s-9} \int }
   +\frac32\partial^2\left(v_x\right)^2\left(t_m+\Dt/2\right)
   \left(t-t_m\right) 
   +\frac{(t-t_m)^2}{2}\int_0^1\Psi\left(\sigma(t-t_m)+t_m\right)d\sigma
   \Bigr)\Bigr]^2\,dx
   \\
   &= \sum_{j=0}^{s-9}\int
   \Bigl[\partial^j \Bigl(
   -\frac32\partial^2\left(v_x\right)^2\left(t_m\right)
   +\frac32\partial^2\left(v_x\right)^2\left(t_m+\Dt/2\right)
   \left(t-t_m\right)
   \\
   &\hphantom{= \sum_{j=0}^{s-9}\int \Bigl(}\quad
   +\frac{(t-t_m)^2}{2}\int_0^1\left(\Phi\left(\sigma(t-t_m)+t_m\right)
     +\Psi\left(\sigma(t-t_m)+t_m\right)\right)d\sigma 
  \Bigr)\Bigr]^2\,dx
   \\
   &\le 2\frac94\sum_{j=0}^{s-9}\int
   \left(\partial^j\left(
       -\partial^2\left(v_x\right)^2\left(t_m\right)
       +\partial^2\left(v_x\right)^2\left(t_m+\Dt/2\right)\right)
     \left(t-t_m\right)\right)^2\, dx
   \\
   &\quad 
   +2\sum_{j=0}^{s-9}\int
   \left(\partial^j \int_0^1\left(\Phi\left(\sigma(t-t_m)+t_m\right)
       +\Psi\left(\sigma(t-t_m)+t_m\right)\right)d\sigma\left(t-t_m\right)^2
   \right)^2\, dx\,
   \\
   &=\frac92 \norm{\Gamma(t)}_{H^{s-9}}^2 
   \left(t-t_m\right)^2 \\
   &\quad
   +2\sum_{j=0}^{s-9}\int
   \left(\partial^j \int_0^1\left(\Phi\left(\sigma(t-t_m)+t_m\right)
       +\Psi\left(\sigma(t-t_m)+t_m\right)\right)d\sigma\right)^2\,dx\,
   \left(t-t_m\right)^4
   \\
   &\le K_\alpha\Dt^2
   \left(t-t_m\right)^2 \\
   &\quad+4\int_0^1\left(\norm{\Phi\left(\sigma(t-t_m)+t_m\right)}_{H^{s-9}}^2
     +\norm{\Psi\left(\sigma(t-t_m)+t_m\right)}_{H^{s-9}}^2\right)d\sigma
   \left(t-t_m\right)^4\\
   &\le K_\alpha\Dt^2\left(t-t_m\right)^2
   +K_\alpha \left(t-t_m\right)^4 \\
   &\le K_\alpha\Dt^4.
 \end{align*}
\end{proof}
Combining the above results, we find that for $t\in [t_{n-1},t_n)$,
\begin{align*}
  \int_0^{t}  &  \norm{H(\sigma)+\widetilde H(\sigma)}_{H^{s-9}}
  \,d\sigma \le \int_0^{t_n} \norm{H(\sigma)+\widetilde H(\sigma)}_{H^{s-9}}
  \,d\sigma 
  \\
  &= \sum_{m=0}^{n-1} \int_{t_m}^{t_{m+1}} \biggl(\, \sum_{j=0}^{s-9}
  \int \left(\partial^j\left(F(\sigma)+G(\sigma)
    \right)+\partial^j\big(F(\sigma+\frac{\Dt}{2})+G(\sigma+\frac{\Dt}{2})
    \big)\right)^2 \, dx \biggr)^{1/2}\,d\sigma 
  \\
  &= \sum_{m=0}^{n-1} \int_{t_m}^{t_{m+1}} \biggl(\, \sum_{j=0}^{s-9}
  \int \left(\partial^j\big(F(\sigma)+G(\sigma+\frac{\Dt}{2})
    \big)+\partial^j\big(F(\sigma+\frac{\Dt}{2})+G(\sigma)
    \big)\right)^2 \, dx \biggr)^{1/2}\,d\sigma 
  \\
  &\le \sqrt{2}\sum_{m=0}^{n-1} \int_{t_m}^{t_{m+1}}
  \biggl(\, \sum_{j=0}^{s-9} \int
  \biggl(\Big(\partial^j\big(F(\sigma)+G(\sigma+\frac{\Dt}{2})\big)\Big)^2
  \\ 
  &\qquad\qquad\qquad\qquad
  +\Big(\partial^j\big(F(\sigma+\frac{\Dt}{2})+G(\sigma)
  \big)\Big)^2 \, dx \biggr)^{1/2}\,d\sigma 
  \\
  &\le \sqrt{2}\sum_{m=0}^{n-1} \int_{t_m}^{t_{m+1}} \biggl[ \biggl(\,
  \sum_{j=0}^{s-9} \int
  \Bigl(\partial^j\big(F(\sigma)+G(\sigma+\frac{\Dt}{2})
  \big)\Bigr)^2 \, dx \biggr)^{1/2} 
  \\
  &\qquad\qquad\qquad\qquad + \biggl(\, \sum_{j=0}^{s-9} \int
  \Big(\partial^j\big(F(\sigma+\frac{\Dt}{2})+G(\sigma) \big)\Bigr)^2
  \, dx \biggr)^{1/2}
  \biggr] \,d\sigma  
  \\
  &= \sqrt{2}\sum_m  
  \int_{t_m}^{t_{m+1/2}} \biggl(\,\sum_{j=0}^{s-9}\int
  \left(\partial^j \left(F(\sigma)+G(\sigma+\Dt/2)\right)\right)^2\,dx
  \biggr)^{1/2} \,d\sigma 
  \\
  &\hphantom{=\sum_m\sqrt{2}}\quad + \int_{t_{m+1}}^{t_{m+3/2}}
  \biggl(\,\sum_{j=0}^{s-9}\int
  \left(\partial^j \left(F(\sigma)+G(\sigma-\Dt/2)\right)\right)^2\,dx
  \biggr)^{1/2} \,d\sigma 
  \\ 
  &\le K_\alpha  \Dt^2 \left(\sum_m \int_{t_m}^{t_{m+1/2}} \,d\sigma +
  \sum_m \int_{t_m+1}^{t_{m+3/2}} \,d\sigma  \right)
  \\
  &\le K_\alpha \Dt^2,
\end{align*}
where we have used Lemma~\ref{lem:cancellation}. This finishes the
estimate for the forcing term.

Next we turn to the estimate of the term $\int_0^t \norm{\tilde
  w(\sigma)(\tilde u(\sigma)-u(\sigma))-w(\sigma)\tilde
  w(\sigma)}_{H^{s-9}}\,d\sigma$ in \eqref{eq:EStrang}.  Here we can
use the estimates from the Godunov splitting to infer that
\begin{equation*}
\norm{w(\sigma)}_{H^{s-9}}+\norm{\tilde w(\sigma)}_{H^{s-9}} \le K_\alpha\Dt, \quad \sigma\in[0,T].
\end{equation*}
From the KdV equation we infer immediately
\begin{equation*}
\norm{\tilde u(\sigma)-u(\sigma)}_{H^{s-9}} 
\le\int_0^{\Dt/2}\norm{(u u_x-u_{xxx})(\sigma+\tau)}_{H^{s-9}}  d\tau  \le K\Dt, \quad \sigma\in[0,T].
\end{equation*}
Thus
\begin{align}
\int_0^t &
\norm{\tilde w(\sigma)(\tilde u(\sigma)-u(\sigma))-w(\sigma)\tilde w(\sigma)}_{H^{s-9}}\,d\sigma\notag  \\
&\le K\int_0^t\Big(\norm{\tilde w(\sigma)}_{H^{s-9}}\norm{(\tilde u-u)(\sigma)}_{H^{s-9}} +
\norm{w(\sigma)}_{H^{s-9}}\norm{\tilde w(\sigma)}_{H^{s-9}} \Big)\,d\sigma \notag\\
&\le K_\alpha \Dt^2.\label{eq:Strang11}
\end{align}

The last term to estimate in  \eqref{eq:EStrang} is $E(0)=w(0)+\tilde
w(0)=w(\Dt/2)$. For $t\le \Dt/2$, we find that
\begin{align*}
  w(t)&=v(t,t)-u(t)=t \int_0^1\big(B(v(st,0))+A(v(t,st))-C(u(st)) \big)ds\\
  &=t \int_0^1\Big(B(u(st))+\int_0^1
  dB(u(st)+\sigma(v(st,0)-u(st)))[v(st,0)-u(st)]d\sigma\\ 
  &\qquad +A(u(st))+\int_0^1
  dA(u(st)+\sigma(v(t,st)-u(st)))[v(t,st)-u(st)]d\sigma\\ 
  &\qquad -(A+B)(u(st)) \Big) ds\\
  &=t \int_0^1 \int_0^1\Big(dB(u(st)+\sigma(v(st,0)-u(st)))[v(st,0)-u(st)]\\
  &\qquad+dA(u(st)+\sigma(v(t,st)-u(st)))[v(t,st)-u(st)] \Big) ds\,d\sigma\\
  &=t \int_0^1
  \int_0^1\Big(dB(u(st)+\sigma(v(st,0)-u(st)))\Bigl[\int_0^{st}
  \frac{d}{d\tau}(v(\tau,0)-u(\tau))\,d\tau\Bigr]
  \\ 
  &\qquad+dA(u(st)+\sigma(v(t,st)-u(st)))\Bigl[\int_0^{st}\frac{d}{d\tau}
  (v(t,\tau)-u(\tau))d\tau+ 
  \int_0^t  \frac{d}{d\tau}v(\tau,0)\,d\tau\Bigr] \Big)\, ds\,d\sigma
  \\
  &=t \int_0^1 \int_0^1\biggl[\int_0^{st}
  \Big(dB(u(st)+\sigma(v(st,0)-u(st)))[B(v(\tau,0))-(A+B)(u(\tau))]
  \\ 
  &\qquad+dA(u(st)+\sigma(v(t,st)-u(st)))[B(v(t,\tau))-(A+B)(u(\tau))]
  \Big)\, d\tau\\ 
  &\qquad +\int_0^t dA(u(st)+\sigma(v(t,st)-u(st)))[B(v(\tau,0))]d\tau
  \biggr] ds\,d\sigma 
  \\
  &=t \int_0^1 \int_0^1\biggl[
  \int_0^{st}
  \Big(\big(u(st)+\sigma(v(st,0)-u(st))\big)
  \big(v(\tau,0)v(\tau,0)_x-u(\tau)u(\tau)_x+u(\tau)_{xxx} 
  \big)\Big)_x  
  \\ 
  &\hphantom{=t \int_0^1 \int_0^1\biggl[\int_0^{st}}\quad
  +\Big(v(t,\tau)v(t,\tau)_x-u(\tau)u(\tau)_x+u(\tau)_{xxx}\Big)_{xxx}\Big)
  d\tau
  \\
  &\hphantom{=t \int_0^1 \int_0^1\Bigl[}\quad 
  +\int_0^t \Big(v(\tau,0)v(\tau,0)_x
  \Big)_{xxx}\,d\tau\biggr]\,ds\,d\sigma.
\end{align*}
Taking the $H^{s-9}$ norm above, using the triangle inequality, the
bounds on $v$ and $u$ and the Cauchy--Schwarz inequality, we find that
each of the above integrands are bounded by $K_\alpha$. Thus
\begin{equation}\label{eq:Strang12}
\begin{aligned}
\norm{w(t)}_{H^{s-9}}&\le t \int_0^1 \int_0^1\Big[\int_0^{st} K_\alpha d\tau+ \int_0^t K_\alpha d\tau \Big]ds\,d\sigma  \\
&\le K_\alpha t^2.
\end{aligned}
\end{equation}
Hence we infer that
\begin{equation} \label{eq:Strang13}
E(0)=\norm{w\left(\Dt/2\right)}_{H^{s-9}}\le  K_\alpha \Dt^2.
\end{equation}
Collecting the estimates from \eqref{eq:Strang11}, \eqref{eq:Strang12}, and \eqref{eq:Strang13}, we find that \eqref{eq:EStrang} reads
\begin{equation} \label{eq:Strang14}
E(t)=\norm{z(t)}_{H^{s-9}}\le  K_\alpha \Dt^2.
\end{equation}
By the triangle inequality
\begin{equation} \label{eq:Strang15}
  \norm{w(t)}_{H^{s-9}}\le 2\norm{w(t)}_{H^{s-9}}\le
  \norm{z(t)}_{H^{s-9}}+\norm{w(t)-\tilde w(t)}_{H^{s-9}}. 
\end{equation}
To estimate the last term on the right-hand side we write:
\begin{equation} \label{eq:last}
\begin{aligned}
\norm{\tilde w(t)- w(t)}_{H^{s-9}}&\le
\norm{w(t_n)- w(t_{n+1/2})}_{H^{s-9}}+\norm{w(t)- w(t_n)}_{H^{s-9}}\\
&\qquad\qquad\qquad\qquad 
 +\norm{w(t+\Dt/2)- w(t_{n+1/2})}_{H^{s-9}}, 
\end{aligned}
\end{equation}
for $t\in [t_n,t_{n+1/2}]$. (Similar expressions hold when $t\in
[t_{n+1/2},t_{n+1}]$.)  We note that
\begin{align}
w(t_{n+1/2})&=
\Phi_A(\frac{\Dt}{2})\Phi_B(\frac{\Dt}{2})v(t_n,t_n)-\Phi_C(\frac{\Dt}{2})u(t_n),\notag 
\\ 
w(t)&=\Phi_A(t-t_n)\Phi_B(t-t_n)v(t_n,t_n)-\Phi_C(t-t_n)u(t_n),\label{eq:last1}\\
w(t+\frac{\Dt}{2})&=\Phi_A(t-t_n)\Phi_B(t-t_n)v(t_{n+1/2},t_{n+1/2})
-\Phi_C(\frac{\Dt}{2})u(t_{n+1/2}), \notag
\end{align}
when $t\in [t_n,t_{n+1/2}]$.  Each of the expressions on the
right-hand side of \eqref{eq:last} needs to be estimated:
\begin{align}
  w(t_{n+1/2})- w(t_n)&=
  \Phi_A\big(\frac{\Dt}{2};\Phi_B(\frac{\Dt}{2};v(t_n,t_n)\big)\big)
  -\Phi_C\big(\frac{\Dt}{2};u(t_n)\big)-(v(t_n,t_n)-u(t_n))
  \notag\\ 
  &=\Phi_A\big(\frac{\Dt}{2};\Phi_B(\frac{\Dt}{2};v(t_n,t_n)\big)\big)
  -\Phi_C\big(\frac{\Dt}{2};v(t_n,t_n)\big)
  \notag\\
  &\qquad 
  +\Phi_C\big(\frac{\Dt}{2};v(t_n,t_n)\big)-\Phi_C\big(\frac{\Dt}{2};u(t_n)\big)
  -(v(t_n,t_n)-u(t_n)) 
  \notag\\
  &=\left(\Phi_A(\frac{\Dt}{2};\Phi_B(\frac{\Dt}{2};v(t_n,t_n)))
  -\Phi_C(\frac{\Dt}{2};v(t_n,t_n))\right)
  \notag\\
  &\qquad 
  +\Bigl(\Phi_C(\frac{\Dt}{2};\dott)-I\Bigr)\circ
  v(t_n,t_n)-\Bigl(\Phi_C(\frac{\Dt}{2};\dott)-I\Bigr)\circ
  u(t_n). \label{eq:nextterm} 
\end{align}
First we find that
\begin{equation*}
\norm{\Phi_A(\frac{\Dt}{2};\Phi_B(\frac{\Dt}{2};v(t_n,t_n)))
  -\Phi_C(\frac{\Dt}{2};v(t_n,t_n))}_{H^{s-9}}
 \le K_\alpha \Dt^2
\end{equation*}
by using \eqref{eq:Strang12} (with $v(t_n,t_n)$ as initial data).

Introduce the function $V=V(x,t)$ satisfying
\begin{equation*}
V_t=V V_x-V_{xxx}, \quad V|_{t=t_n}=v(t_n,t_n).
\end{equation*}
Then the very last line of \eqref{eq:nextterm} can be written as
\begin{equation*}
\begin{aligned}
&\Bigl(\Phi_C(\frac{\Dt}{2};\dott)-I\Bigr)\circ
  v(t_n,t_n)-\Bigl(\Phi_C(\frac{\Dt}{2};\dott)-I\Bigr)\circ u(t_n) 
 \\
&\quad =\int_{t_n}^{t_{n+1/2}}\big(V_t(\sigma)- u_t(\sigma)\big)d\sigma\\
&\quad = \int_{t_n}^{t_{n+1/2}}\Big(\frac12 (V(\sigma)^2)_x-\frac12 (u(\sigma)^2)_x-V_{xxx}(\sigma)+u_{xxx}(\sigma) \Big)d\sigma\\ 
&\quad = \int_{t_n}^{t_{n+1/2}}\Big(\frac12 \big((V(\sigma)+u(\sigma))(V(\sigma)-u(\sigma))\big)_x-(V(\sigma)-u(\sigma))_{xxx}\Big)d\sigma.
\end{aligned}
\end{equation*}
Taking the $H^{s-9}$ norm we find
\begin{equation*}
\begin{aligned}
&\norm{\Bigl(\Phi_C(\frac{\Dt}{2};\dott)-I\Bigr)\circ
  v(t_n,t_n)-\Bigl(\Phi_C(\frac{\Dt}{2};\dott)-I\Bigr)\circ u(t_n)}_{H^{s-9}}
\\
&\quad \le \int_{t_n}^{t_{n+1/2}}\Big(\frac12 \norm{(V(\sigma)+u(\sigma))(V(\sigma)-u(\sigma))}_{H^{s-8}}+\norm{V(\sigma)-u(\sigma)}_{H^{s-6}}\Big)d\sigma\\
&\quad \le \int_{t_n}^{t_{n+1/2}}\Big(K \big(\norm{V(\sigma)}_{H^{s-8}}+ \norm{u(\sigma)}_{H^{s-8}}\big)\norm{V(\sigma)-u(\sigma)}_{H^{s-8}}+\norm{V(\sigma)-u(\sigma)}_{H^{s-6}}\Big)d\sigma\\
&\quad \le K_\alpha  \int_{t_n}^{t_{n+1/2}}\norm{V(\sigma)-u(\sigma)}_{H^{s-6}}d\sigma.
\end{aligned}
\end{equation*}
By the $H^k$ stability of the KdV equation 
\begin{equation*}
\norm{V(\sigma)-u(\sigma)}_{H^{s-6}}\le
K\norm{v(t_n,t_n)-u(t_n)}_{H^{s-6}}\le K_\alpha\Dt, 
\end{equation*}
since by the arguments of Section~\ref{sec:kdv_split} we have the
estimate $\norm{w(t_n)}_{H^{s-6}}\le \norm{w(t_n)}_{H^{s-3}}\le
K_\alpha \Dt$.
Therefore 
\begin{equation*}
\norm{\Bigl(\Phi_C(\frac{\Dt}{2};\dott)-I\Bigr)\circ
  v(t_n,t_n)-\Bigl(\Phi_C(\frac{\Dt}{2};\dott)-I\Bigr)\circ u(t_n)}_{H^{s-9}}
\le K_\alpha \Dt^2.
\end{equation*}
Thus we have shown that
\begin{equation} \label{eq:Strang15Aa}
  \norm{w(t_{n+1/2})- w(t_n)}_{H^{s-9}}\le K_\alpha\Dt^2.
\end{equation}
The other terms on the right hand side of \eqref{eq:last} can be
estimated in the same manner, using the expressions \eqref{eq:last1}.
Thus we conclude that (cf.~\eqref{eq:Strang15})
\begin{equation} \label{eq:Strang15A}
  \norm{w(t)}_{H^{s-9}}\le K_\alpha\Dt^2.
\end{equation}

If $t,\tau\in [t_n,t_{n+1/2}]$, we have
$\norm{v(t,\tau)}_{H^{\hat{k}}}=\norm{v(t,t)}_{H^{\hat{k}}}$, and  if
$t,\tau\in [t_{n+1/2},t_n]$, then an estimate analogous to
\eqref{eq:burgerhsdev}  shows that 
$$
\abs{\vphantom{\Bigm|}\norm{v(t,\tau)}_{H^{\hat{k}}}-\norm{v(t,t)}_{H^{\hat{k}}}}
\le K_\alpha 
\abs{t-\tau}.
$$
The rest of the argument follows the procedure for the Godunov
splitting.  Now  $s-9\ge \hat{k}$ is the same as $s\ge
17$. Assuming this, we get
\begin{align*}
 \norm{v(t,\tau)}_{H^{\hat{k}}} &\le \norm{v(t,t)}_{H^{\hat{k}}} + 
 \abs{\vphantom{\Bigm|}
   \norm{v(t,\tau)}_{H^{\hat{k}}}-\norm{v(t,t)}_{H^{\hat{k}}}}\\ 
 &\le K + K_\alpha\Dt^2 + K_\alpha\Dt.
\end{align*}
Choosing $\alpha$ such that $K\le \alpha/4$, and then $\Dt$ such that $K_\alpha
\Dt(\Dt + 1)\le \alpha/4$  implies that $\norm{v(t,\tau)}_{H^{\hat{k}}}
\le \alpha/2$. Hence by the bootstrap lemma and \eqref{eq:Strang15}, 
the following holds:
\begin{theorem}
  \label{thm:Strangsplit} Fix $T>0$.  Let $u_0\in H^s$ for some $s\ge
  17$.  Then for $\Dt$ sufficiently small we have
  $$
  \norm{v(t,t)-u(t)}_{H^{s-9}} \le K \Dt^2, \quad t\in [0,T],
  $$
  where the constant $K$ depends on $u_0$, $s$ and $T$ only.
\end{theorem}

\end{document}